\documentclass[final,1p,times,12pt]{article}

\usepackage{enumerate}
\usepackage{epsfig,alltt}

\usepackage{amssymb}
\usepackage{amsfonts}
\usepackage{graphicx}
\usepackage{amsmath}
\usepackage{amsthm}
\usepackage{bm}
\usepackage{natbib}

\usepackage[utf8]{inputenc}
\usepackage[T1]{fontenc}

\usepackage{makeidx}

\usepackage{url}
\usepackage{float} 

\theoremstyle{remark}

\usepackage[english]{babel}

\usepackage{amsfonts}
\usepackage{amssymb}
\usepackage{geometry}
\usepackage{color}

\newcommand{\begeq}[1]{\begin{equation} \label{#1}}
\newcommand{\fineq}{\end{equation}}

\def\bE{{\Bbb E}}                  % Esperance mathematique
\title{Optimal bandwidth selection for semi-recursive kernel regression estimators}
\author{Yousri Slaoui\\
Universit\'e de Poitiers}

\begin{document}

\newtheorem{theor}{Theorem}
\newtheorem{prop}{Proposition}
\newtheorem{lemma}{Lemma}
\newtheorem{lem}{Lemma}
\newtheorem{coro}{Corollary}

\newtheorem{prof}{Proof}
\newtheorem{defi}{Definition}
\newtheorem{rem}{Remark}

\date{ }
\maketitle

\textit{Abstract}: In this paper we propose an automatic selection of the bandwidth of the semi-recursive kernel estimators of a regression function defined by the stochastic approximation algorithm. We showed that, using the selected bandwidth and some special stepsizes, the proposed semi-recursive estimators will be very competitive to the nonrecursive one in terms of estimation error but much better in terms of computational costs. We corroborated these theoretical results through simulation study and a real dataset.\\

%2010 \textit{Mathematics Subject Classification}: Primary 62E20, 62L20.\\
\textit{Key words and phrases:} Nonparametric regression; Stochastic approximation algorithm; Smoothing, curve fitting.\\
{\bf Mathematics Subject Classification (2010)} Primary $62$G$08$. $62$L$20$. $65$D$10$.
\section{Introduction}
In recent years, there has been a lot of interest in big data. In such a large sample data context, building a semi-recursive estimator which does not require to store all the data in memory and can be updated easily in order to deal with online data is of great interest.\\
In the framework of the nonparametric kernel estimators, the bandwidth selection methods studied in the literature can be divided into three broad classes: the cross-validation techniques, the plug-in ideas and the bootstrap. A detailed comparison of the three practical bandwidth selection can be found in \citet{Del04}. They concluded that chosen appropriately plug-in and bootstrap selectors both outperform the cross-validation bandwidth, and that neither of the two can be claimed to be better in all cases. Recently, plug-in bandwidth selection method for recursive kernel density estimators defined by stochastic approximation method have been done by~\citet{Sla14a} and for recursive kernel distribution estimators have been done by~\citet{Sla14b}. In this paper, we developed a specific plug-in bandwidth selection method of the semi-recursive kernel estimators of a regression function defined by stochastic approximation method.

Let $\left(X_1,Y_1\right),\ldots ,\left(X_n,Y_n\right)$ be independent, identically distributed pairs of random variables with joint density function $g\left(x,y\right)$, and let $f$ denote the probability density of $X$. In order to construct a stochastic
algorithm for the estimation of the regression function $a:x\mapsto
\bE(Y|X=x)f\left(x\right)$ at a point $x$, we define an algorithm of search of the zero of the function $h : y\to a(x)-y$. Following Robbins-Monro's procedure, this algorithm is defined by setting $a_0(x)\in \mathbb{R}$, and, for all $n\geq 1$, 
\begin{eqnarray*}
a_n\left(x\right)=a_{n-1}\left(x\right)+\beta_nW_n,
\end{eqnarray*}
where $W_n(x)$ is an "observation" of the function $h$ at the point $a_{n-1}(x)$, and the stepsize $\left(\beta_n\right)$ is a sequence of positive real numbers that goes to zero. To define $W_n(x)$, we follow the approach of \citet{Rev73,Rev77},~\citet{Tsy90} and of \citet{Mok09a,Mok09b} and introduces a kernel $K$ (that is, a function satisfying $\int_{\mathbb{R}} K(x)dx=1$), and a bandwidth $(h_n)$ (that is,
a sequence of positive real numbers that goes to zero), and sets $W_n(x)=h_n^{-1}Y_nK\left(h_n^{-1}\left(x-X_n\right)\right)-a_{n-1}(x)$.
Then, the estimator $a_n$ to recursively estimate the function $a$ at the point $x$ can be written as
\begin{eqnarray}\label{eq:an}
a_n\left(x\right)=\left(1-\beta_n\right)a_{n-1}\left(x\right)+\beta_nh_n^{-1}Y_nK\left(h_n^{-1}\left[x-X_n\right]\right).
\end{eqnarray}
This estimator was proposed by~\citet{Sla15c} to estimate recursively the regression function with a fixed design setting. The recursive property~(\ref{eq:an}) is particularly useful in large sample size since $a_n$ can be easily updated with each additional observation.\\
Let us underline that, we consider $a_0\left(x\right)=0$ and we let $Q_n=\prod_{j=1}^n\left(1-\beta_j\right)$, then it follows from~(\ref{eq:an}) that, one can estimate $a$ recursively at the point $x$ by 
\begin{eqnarray*}
a_n\left(x\right)
&=&Q_n\sum_{k=1}^nQ_k^{-1}\beta_kh_k^{-1}Y_kK\left(\frac{x-X_k}{h_k}\right).
\end{eqnarray*} 
Moreover, we use the estimator introduced in~\citet{Mok09a} to estimate recursively the density $f$ at the point $x$
\begin{eqnarray}\label{eq:fn}
f_n\left(x\right)=\left(1-\gamma_n\right)f_{n-1}\left(x\right)+\gamma_nh_n^{-1}K\left(h_n^{-1}\left[x-X_n\right]\right),
\end{eqnarray}
where the stepsize $\left(\gamma_n\right)$ is a sequence of positive real numbers that goes to zero. Let us underline that we consider  $f_0\left(x\right)=0$, and we let $\Pi_n=\prod_{j=1}^n\left(1-\gamma_j\right)$, then it follows from~(\ref{eq:fn}) that, one can estimate $f$ recursively at the point $x$ by 
\begin{eqnarray*}
f_n\left(x\right)
&=&\Pi_n\sum_{k=1}^n\Pi_k^{-1}\gamma_kh_k^{-1}K\left(\frac{x-X_k}{h_k}\right).
\end{eqnarray*}
Then, we consider the semi-recursive estimator for the regression function $r$ at the point $x$
\begin{eqnarray}\label{eq:rn}
r_n\left(x\right)=\left\{\begin{array}{cc}
\frac{a_n\left(x\right)}{f_n\left(x\right)} & \mbox{if}\quad f_n\left(x\right)\not =0,\\
0  &\mbox{otherwise}.
\end{array}\right.
\end{eqnarray} 
Moreover, we show that the optimal bandwidth which minimize the $\mathbb{E}\int_{\mathbb{R}}\left[r_n\left(x\right)-r\left(x\right)\right]^2dx$ of $r_n$ depends on the choice of the stepsizes $\left(\gamma_n\right)$ and $\left(\beta_n\right)$; we show in particular that under some conditions of regularity of $r$ and using the stepsizes $\left(\gamma_n,\beta_n\right)=\left(n^{-1},n^{-1}\right)$, the bandwidth $\left(h_n\right)$ must equal 
\begin{eqnarray*}
\left(\left(\frac{3}{10}\right)^{1/5}\left\{\frac{\int_{\mathbb{R}}Var\left[Y^2\vert X=x\right]f^{-1}\left(x\right)dx}{\int_{\mathbb{R}}\left(a^{\left(2\right)}\left(x\right)-r\left(x\right)f^{\left(2\right)}\left(x\right)\right)^2f^{-2}\left(x\right)dx}\right\}^{1/5}\left\{\frac{\int_{\mathbb{R}}K^2\left(z\right)dz}{\left(\int_{\mathbb{R}}z^2K\left(z\right)dz\right)^2}\right\}^{1/5}n^{-1/5}\right).
\end{eqnarray*}
The first aim of this paper is to propose an automatic selection of such bandwidth through a plug-in method, and the second aim is to give the conditions under which the semi-recursive estimator $r_n$ will be approximately similar to the nonrecursive kernel regression estimators introduced by~\citet{Nad64} and~\citet{Wat64}, and defined as
\begin{eqnarray}\label{eq:Nad}
\widetilde{r}_n\left(x\right)=\left\{\begin{array}{cc}
\frac{\widetilde{a}_n\left(x\right)}{\widetilde{f}_n\left(x\right)} & \mbox{if}\quad \widetilde{f}_n\left(x\right)\not =0\\
0  &\mbox{otherwise},
\end{array}\right.
\end{eqnarray}
with 
\begin{eqnarray*}
\widetilde{a}_n\left(x\right)=\frac{1}{nh_n}\sum_{i=1}^nY_iK\left(\frac{x-X_i}{h_n}\right)\quad \mbox{and}\quad \widetilde{f}_n\left(x\right)=\frac{1}{nh_n}\sum_{i=1}^nK\left(\frac{x-X_i}{h_n}\right).
\end{eqnarray*}
The applications results given in Section~\ref{section:app} corroborate these theoretical results. 
The remainder of the paper is organized as follows. In Section~\ref{section:2}, we state our main results. Section~\ref{section:app} is devoted to our application results, first by simulation (subsection~\ref{subsection:simu}) and second using a real dataset (subsection~\ref{subsection:real}). We conclude the article in Section~\ref{section:conclusion}. Appendix~\ref{section:proofs} gives the proof of our theoretical results. 
\section{Assumptions and main results} \label{section:2}
We define the following class of regularly varying sequences.
\begin{defi}
Let $\gamma \in \mathbb{R} $ and $\left(v_n\right)_{n\geq 1}$ be a nonrandom positive sequence. We say that $\left(v_n\right) \in \mathcal{GS}\left(\gamma \right)$ if
\begin{eqnarray}\label{eq:5}
\lim_{n \to +\infty} n\left[1-\frac{v_{n-1}}{v_{n}}\right]=\gamma .
\end{eqnarray}
\end{defi}
Condition~(\ref{eq:5}) was introduced by \citet{Gal73} to define regularly varying sequences (see also \citet{Boj73}) and by \citet{Mok07} in the context of stochastic approximation algorithms. Noting that the acronym $\mathcal{GS}$ stand for (Galambos and Seneta). Typical sequences in $\mathcal{GS}\left(\gamma \right)$ are, for $b\in \mathbb{R}$, $n^{\gamma}\left(\log n\right)^{b}$, $n^{\gamma}\left(\log \log n\right)^{b}$, and so on. \\
In this section, we investigate asymptotic properties of the proposed estimators~(\ref{eq:rn}). The assumptions to which we shall refer are the following:

\begin{description}
\item(A1) $K:\mathbb{R}\rightarrow \mathbb{R}$ is a continuous, bounded function satisfying $\int_{\mathbb{R}}K\left( z\right) dz=1$, and, $\int_{\mathbb{R}}zK\left(z\right)=0$ and $\int_{\mathbb{R}}z^2K\left(z\right)<\infty$. 
\item(A2) $i)$ $\left(\beta_n\right)\in \mathcal{GS}\left(-\beta\right)$ with $\beta\in \left]1/2,1\right]$. \\
  $ii)$ $\left(h_n\right)\in \mathcal{GS}\left(-a\right)$ with $a\in \left]0,1\right[$.\\
  $iii)$ $\lim_{n\to \infty} \left(n\beta_n\right)\in \left]\min\left\{2a,\left(\beta-a\right)/2\right\},\infty\right]$.
\item(A3) $i)$ $g\left(s,t\right)$ is twice continuously differentiable with
respect to $s$.\\
 $ii)$ For $q\in\left\{0,1,2\right\}$, $s \mapsto
\int_{\mathbb{R}}t^qg\left(s,t\right)dt$ is a bounded function continuous at
$s=x$.\\ 
 For $q\in \left[2,3\right]$, $s \mapsto
\int_{\mathbb{R}}\left|t\right|^qg\left(s,t\right)dt$ is a bounded function.\\
$iii)$ For $q\in \left\{0,1\right\}$,
$\int_{\mathbb{R}}\left|t\right|^q\left|\frac{\partial g}{\partial
x}\left(x,t\right)\right|dt<\infty$, and $s\mapsto
\int_{\mathbb{R}}t^q\frac{\partial^2 g}{\partial s^2}\left(s,t\right)dt$ is a
bounded function continuous at $s=x$. 
\end{description}
Assumption $\left(A2\right)(iii)$ on the limit of $\left(n\beta_n\right)$ as $n$ goes to infinity is standard in the framework of stochastic approximation algorithms. It implies in particular that the limit of $\left(\left[n\beta_n\right]^{-1}\right)$ is finite. For simplicity, we introduce the following notations:
\begin{eqnarray}
\xi&=&\lim_{n\to \infty}\left(n\beta_n\right)^{-1},\label{eq:xi}
\end{eqnarray}
\begin{eqnarray*}
R\left(K\right)=\int_{\mathbb{R}}K^2\left(z\right)dz,\quad \quad\mu_j\left(K\right)=\int_{\mathbb{R}}z^jK\left(z\right)dz,\\
\Theta\left(K\right)=R\left(K\right)^{4/5}\mu_2\left(K\right)^{2/5},\quad
I_1=\int_{\mathbb{R}}\left(a^{\left(2\right)}\left(x\right)\right)^2f\left(x\right)dx,\\
I_2=\int_{\mathbb{R}}a^{\left(2\right)}\left(x\right)f^{\left(2\right)}\left(x\right)r\left(x\right)f\left(x\right)dx,\quad
I_3=\int_{\mathbb{R}}\left(f^{\left(2\right)}\left(x\right)\right)^2r^2\left(x\right)f\left(x\right)dx,\\
I_4=\int_{\mathbb{R}}\mathbb{E}\left[Y^2\vert X=x\right]f^2\left(x\right)dx,\quad
I_5=\int_{\mathbb{R}}r^2\left(x\right)f^2\left(x\right)dx,
\end{eqnarray*} 
where $L^{\left(2\right)}\left(x\right)$ is the second derivative of the function $L$ at a point $x$.
In this section, we explicit the choice of $\left(h_n\right)$ through a plug-in method, which minimize the Mean Weighted Integrated Squared Error $MWISE$ of the semi-recursive estimators~(\ref{eq:rn}), in order to provide a comparison with the nonrecursive estimator~(\ref{eq:Nad}). Moreover, it was shown in~\citet{Mok09a} and considered in~\citet{Sla13} that to minimize the Mean Integrated Squared Error $MISE$ of $f_n$ ($MISE\left[f_n\right]=\mathbb{E}\int_{\mathbb{R}}\left[f_n\left(x\right)-f\left(x\right)\right]^2dx$), the stepsize $\left(\gamma_n\right)$ must be chosen in $\mathcal{GS}\left(-1\right)$ and must satisfy $\lim_{n\to\infty}n\gamma_n=1$. We consider here the case $\left(\gamma_n\right)=\left(n^{-1}\right)$. Our first result is the following proposition, which gives the bias and the variance of $r_n$ in the special case of $\left(\gamma_n\right)=\left(n^{-1}\right)$.
\begin{prop}[Bias and variance of $r_n$]\label{prop:bias:var:rn}
Let Assumptions $\left(A1\right)-\left(A3\right)$ hold, and suppose that the stepsize $\left(\gamma_n\right)=\left(n^{-1}\right)$
\begin{enumerate}
\item If $a\in ]0, \beta/5]$, then
\begin{eqnarray}\label{bias:rn}
\mathbb{E}\left[r_n\left(x\right)\right]-r\left(x\right)=\frac{1}{2f\left(x\right)}\left(\frac{a^{\left(2\right)}\left(x\right)}{\left(1-2a\xi\right)}-\frac{r\left(x\right)f^{\left(2\right)}\left(x\right)}{\left(1-2a\right)}\right)h_n^2\mu_2\left(K\right)+o\left(h_n^2\right).
\end{eqnarray}
If $a\in  ]\beta/5, 1[$, then
\begin{eqnarray}\label{bias:rn:bis}
\mathbb{E}\left[r_n\left(x\right)\right]-r\left(x\right)=o\left(\sqrt{\beta_nh_n^{-1}}\right).
\end{eqnarray}
\item If $a\in [\beta/5, 1[$, then
\begin{eqnarray}
Var\left[r_n\left(x\right)\right]&=&\frac{\beta_n}{h_n}\left\{\frac{\mathbb{E}\left[Y^2\vert X=x\right]}{\left(2-\left(\beta-a\right)\xi\right)f\left(x\right)}-\left(\frac{2\xi}{1+a\xi}-\frac{\xi}{1+a}\right)\frac{r^2\left(x\right)}{f\left(x\right)}\right\}\nonumber\\
&&R\left(K\right)+o\left(\frac{\beta_n}{h_n}\right).\label{var:rn}
\end{eqnarray}
If $a\in ]0,\beta/5[$, then
\begin{eqnarray}\label{var:rn:rep}
Var\left[r_n\left(x\right)\right]=o\left(h_n^4\right).
\end{eqnarray}
\item If $\lim_{n\to \infty}\left(n\beta_n\right)>\max\left\{2a, \left(a-\beta\right)/2\right\}$, then~(\ref{bias:rn}) and~(\ref{var:rn}) hold simultaneously.
\end{enumerate}
\end{prop}
The bias and the variance of the estimator $r_n$ defined by the stochastic approximation algorithm~(\ref{eq:rn}) then heavily depend on the choice of the stepsizes $\left(\gamma_n\right)$ and $\left(\beta_n\right)$.\\
Let us first state the following theorem, which gives the weak convergence rate of the estimator $r_n$ defined in~(\ref{eq:rn}) in the case of $\left(\gamma_n\right)=\left(n^{-1}\right)$.
\begin{theor}[Weak pointwise convergence rate]\label{theo:TLC1}
Let Assumptions $\left(A1\right)-\left(A3\right)$ hold, and suppose that $\left(\gamma_n\right)=\left(n^{-1}\right)$.
\begin{enumerate}
\item If there exists $c\geq 0$ such that $\beta_n^{-1}h_n^5\to c$, then
\begin{eqnarray*}
\sqrt{\beta_n^{-1}h_n}\left(r_{n}\left(x\right)-r\left( x\right) \right) \stackrel{\mathcal{D}}{\rightarrow} &
\mathcal{N}\left(\sqrt{c}B_{a,\xi}^{\left(1\right)}\left(x\right),V_{a,\xi,\beta}^{\left(1\right)}\left(x\right)\right),
\end{eqnarray*}
where
\begin{eqnarray*}
B_{a,\xi}^{\left(1\right)}\left(x\right)&=&\frac{1}{2f\left(x\right)}\left(\frac{a^{\left(2\right)}\left(x\right)}{\left(1-2a\xi\right)}-\frac{r\left(x\right)f^{\left(2\right)}\left(x\right)}{\left(1-2a\right)}\right)\mu_2\left(K\right),\\
V_{a,\xi,\beta}^{\left(1\right)}\left(x\right)&=&\left\{\frac{\mathbb{E}\left[Y^2\vert X=x\right]}{\left(2-\left(\beta-a\right)\xi\right)f\left(x\right)}-\left(\frac{2\xi}{1+a\xi}-\frac{\xi}{1+a}\right)\frac{r^2\left(x\right)}{f\left(x\right)}\right\}R\left(K\right).
\end{eqnarray*}
\item If $nh_{n}^{5} \rightarrow \infty $, then  
\begin{eqnarray*}
\frac{1}{h_{n}^{2}}\left(r_{n}\left(x\right)-r\left(x\right) \right) \stackrel{\mathbb{P}}{\rightarrow } B_{a,\xi}^{\left(1\right)}\left(x\right),
\end{eqnarray*}
\end{enumerate}
where $\stackrel{\mathcal{D}}{\rightarrow}$ denotes the convergence in distribution, $\mathcal{N}$ the Gaussian-distribution and $\stackrel{\mathbb{P}}{\rightarrow}$ the convergence in probability.
\end{theor}
The following corollary gives the weak convergence rate of $r_n$ in the two special cases; $\left(\gamma_n,\beta_n\right)=\left(n^{-1},n^{-1}\right)$ and $\left(\gamma_n,\beta_n\right)=\left(n^{-1},\left(1-a\right)n^{-1}\right)$ respectively. 
\begin{coro}[Weak pointwise convergence rate]\label{coro:TLC1}
Let Assumptions $\left(A1\right)-\left(A3\right)$ hold. 
\begin{enumerate}
\item If we suppose that the stepsizes $\left(\gamma_n,\beta_n\right)=\left(n^{-1},n^{-1}\right)$ and if there exists $c\geq 0$ such that $nh_n^5\to c$, then
\begin{eqnarray*}
\sqrt{nh_n}\left(r_{n}\left(x\right)-r\left( x\right) \right) \stackrel{\mathcal{D}}{\rightarrow} &
\mathcal{N}\left(\sqrt{c}B_{a,1}^{\left(1\right)}\left(x\right),V_{a,1,1}^{\left(1\right)}\left(x\right)\right).
\end{eqnarray*}
\item If we suppose that the stepsizes    
$\left(\gamma_n,\beta_n\right)=\left(n^{-1},\left(1-a\right)n^{-1}\right)$, and if there exists $c\geq 0$ such that $nh_n^5\to c$, then
\begin{eqnarray*}
\sqrt{nh_n}\left(r_{n}\left(x\right)-r\left( x\right) \right) \stackrel{\mathcal{D}}{\rightarrow} &
\mathcal{N}\left(\sqrt{c}B_{a,\left(1-a\right)^{-1}}^{\left(1\right)}\left(x\right),V_{a,\left(1-a\right)^{-1},1}^{\left(1\right)}\left(x\right)\right).
\end{eqnarray*}
\end{enumerate}
\end{coro} 
In order to measure the quality of our semi-recursive estimator~(\ref{eq:rn}) in the case when the stepsize $\left(\gamma_n\right)$ is chosen to minimize the $MISE$ of $f_n$, we use the following quantity, 
\begin{eqnarray*}
MWISE\left[r_n\right]&=&\mathbb{E}\int_{\mathbb{R}}\left[r_n\left(x\right)-r\left(x\right)\right]^2f^3\left(x\right)dx\nonumber\\
&=&\int_{\mathbb{R}}\left(\mathbb{E}\left(r_n\left(x\right)\right)-r\left(x\right)\right)^2f^3\left(x\right)dx+\int_{\mathbb{R}}Var\left(r_n\left(x\right)\right)f^3\left(x\right)dx.
\end{eqnarray*}
The following proposition gives the $MWISE$ of the semi-recursive estimators defined in~(\ref{eq:rn}) in the case when $\left(\gamma_n\right)$ is chosen to minimize the $MISE$ of $f_n$.
\begin{prop}[$MWISE$ of $r_n$]\label{prop:MISE:rn}
Let Assumptions $\left(A1\right)-\left(A3\right)$ hold, and suppose that $\left(\gamma_n\right)=\left(n^{-1}\right)$.
\begin{enumerate}
\item If $a\in  ]0, \beta/5[$, then
\begin{eqnarray*}
MWISE\left[r_n\right]=\frac{1}{4}\left(\frac{I_1}{\left(1-2a\xi\right)^2}+\frac{I_3}{\left(1-2a\right)^2}-2\frac{I_2}{\left(1-2a\right)\left(1-2a\xi\right)}\right)h_n^4\mu_2^2\left(K\right)+o\left(h_n^4\right).
\end{eqnarray*}
\item If $a=\beta/5$, then
\begin{eqnarray*}
MWISE\left[r_n\right]&=&\frac{\beta_n}{h_n}\left(\frac{I_4}{\left(2-\left(\beta-a\right)\xi\right)}-\left(\frac{2\xi}{1+a\xi}-\frac{\xi}{1+a}\right)I_5\right)R\left(K\right)\\
&&+\frac{1}{4}\left(\frac{I_1}{\left(1-2a\xi\right)^2}+\frac{I_3}{\left(1-2a\right)^2}-2\frac{I_2}{\left(1-2a\right)\left(1-2a\xi\right)}\right)h_n^4\mu_2^2\left(K\right)\\
&&+o\left(h_n^4\right).
\end{eqnarray*}
\item If $a\in ]\beta/5, 1[$, then
\begin{eqnarray*}
MWISE\left[r_n\right]&=&\frac{\beta_n}{h_n}\left(\frac{I_4}{\left(2-\left(\beta-a\right)\xi\right)}-\left(\frac{2\xi}{1+a\xi}-\frac{\xi}{1+a}\right)I_5\right)R\left(K\right)+o\left(\frac{\beta_n}{h_n}\right).
\end{eqnarray*}
\end{enumerate}
\end{prop}
The following corollary indicates that the bandwidth which minimizes the $MWISE$ of $r_n$ depends on the stepsize $\left(\beta_n\right)$ and then the corresponding $MWISE$ depends also on the stepsize $\left(\beta_n\right)$.
\begin{coro}\label{Coro:hn:MISE}
Let Assumptions $\left(A1\right)-\left(A3\right)$ hold, and suppose that $\left(\gamma_n\right)=\left(n^{-1}\right)$. To minimize the $MWISE$ of $r_n$, the stepsize $\left(\beta_n\right)$ must be chosen in $\mathcal{GS}\left(-1\right)$, the bandwidth $\left(h_n\right)$ must equal 
\begin{eqnarray*}
\left(\left\{\frac{\frac{I_4}{\left(2-\left(\beta-a\right)\xi\right)}-\left(\frac{2\xi}{1+a\xi}-\frac{\xi}{1+a}\right)I_5}{\frac{I_1}{\left(1-2a\xi\right)^2}+\frac{I_3}{\left(1-2a\right)^2}-2\frac{I_2}{\left(1-2a\right)\left(1-2a\xi\right)}}\right\}^{1/5}\left\{\frac{R\left(K\right)}{\mu_2^2\left(K\right)}\right\}^{1/5}\beta_n^{1/5}\right).
\end{eqnarray*}
Then, we have
\begin{eqnarray*} MWISE\left[r_n\right]&=&\frac{5}{4}\left(\frac{I_4}{\left(2-\left(\beta-a\right)\xi\right)}-\left(\frac{2\xi}{1+a\xi}-\frac{\xi}{1+a}\right)I_5\right)^{4/5}\\
&&\times\left(\frac{I_1}{\left(1-2a\xi\right)^2}+\frac{I_3}{\left(1-2a\right)^2}-2\frac{I_2}{\left(1-2a\right)\left(1-2a\xi\right)}\right)^{1/5}\Theta\left(K\right)\beta_n^{4/5}\\
&&+o\left(\beta_n^{4/5}\right).
\end{eqnarray*}
\end{coro}
The following corollary shows that, for a special choice of the stepsize $\left(\beta_n\right)=\left(\beta_0n^{-1}\right)$, which fulfilleds that $\lim_{n\to \infty}n\beta_n=\beta_0$ and that $\left(\beta_n\right)\in \mathcal{GS}\left(-1\right)$, the optimal value for $h_n$ depends on $\beta_0$ and then the corresponding $MWISE$ depend on $\beta_0$.
\begin{coro}\label{Coro:hn:MISE:bis}
Let Assumptions $\left(A1\right)-\left(A3\right)$ hold, and suppose that $\left(\gamma_n\right)=\left(n^{-1}\right)$. To minimize the $MWISE$ of $r_n$, the stepsize $\left(\beta_n\right)$ must be chosen in $\mathcal{GS}\left(-1\right)$, $\lim_{n\to \infty}n\beta_n=\beta_0$, the bandwidth $\left(h_n\right)$ must equal 
\begin{eqnarray}\label{hoptim}
\left(\left(\frac{\beta_0-2/5}{2}\right)^{1/5}\left(\frac{I_4-\frac{\left(7\beta_0-1\right)\left(\beta_0-2/5\right)}{3\beta_0^2\left(\beta_0+1/5\right)}I_5}{I_1+\frac{25}{9}\left(\frac{\beta_0-2/5}{\beta_0}\right)^2I_3-\frac{10}{3}\left(\frac{\beta_0-2/5}{\beta_0}\right)I_2}\right)^{1/5}\left\{\frac{R\left(K\right)}{\mu_2^2\left(K\right)}\right\}^{1/5}n^{-1/5}\right),
\end{eqnarray}
and we then have
\begin{eqnarray}\label{MWISE:hoptim}
MWISE\left[r_n\right]&=&\frac{5}{4}\frac{1}{2^{4/5}}\frac{\beta_0^2}{\left(\beta_0-2/5\right)^{6/5}}\left(I_4-\frac{\left(7\beta_0-1\right)\left(\beta_0-2/5\right)}{3\beta_0^2\left(\beta_0+1/5\right)}I_5\right)^{4/5}\nonumber\\
&&\times \left(I_1+\frac{25}{9}\left(\frac{\beta_0-2/5}{\beta_0}\right)^2I_3-\frac{10}{3}\left(\frac{\beta_0-2/5}{\beta_0}\right)I_2\right)^{1/5}\Theta\left(K\right)n^{-4/5}
\nonumber\\
&&+o\left(n^{-4/5}\right).
\end{eqnarray}
Moreover, the minimum of $\beta_0^2\left(\beta_0-2/5\right)^{-6/5}$ is reached at $\beta_0=1$, then the bandwidth $\left(h_n\right)$ must equal 
\begin{eqnarray}\label{hoptim:gamma}
\left(\left(\frac{3}{10}\right)^{1/5}\left(\frac{I_4-I_5}{I_1+I_3-2I_2}\right)^{1/5}\left\{\frac{R\left(K\right)}{\mu_2^2\left(K\right)}\right\}^{1/5}n^{-1/5}\right),
\end{eqnarray}
and we then have
\begin{eqnarray}\label{MISE:gamma:MISE} 
MWISE\left[r_n\right]
&=&\frac{5}{4}\frac{1}{2^{4/5}}\left(\frac{5}{3}\right)^{6/5}\left(I_4-I_5\right)^{4/5}\times \left(I_1+I_3-2I_2\right)^{1/5}\Theta\left(K\right)n^{-4/5}\nonumber\\
&&+o\left(n^{-4/5}\right).
\end{eqnarray}
\end{coro}
In order to estimate the optimal bandwidth (\ref{hoptim:gamma}), we must estimate $I_1$, $I_2$, $I_3$, $I_4$ and $I_5$. We followed the approach of \citet{Alt95} and \citet{Sla14a,Sla14b}, which is called the plug-in estimate, and we use the following kernel estimators of $I_1$, $I_2$, $I_3$, $I_4$ and $I_5$:

\begin{eqnarray}
\widehat{I}_1&=&\frac{Q_n^2}{n}\sum_{\substack{i,j,k=1\\j\not=k}}^nQ_j^{-1}Q_k^{-1}\beta_j\beta_kb_j^{-3}b_k^{-3}K_b^{\left(2\right)}\left(\frac{X_i-X_j}{b_j}\right)K_b^{\left(2\right)}\left(\frac{X_i-X_k}{b_k}\right)Y_jY_k,\label{I1:rec}\\
\widehat{I}_2&=&\frac{\Pi_nQ_n}{n}\sum_{\substack{i,j,k=1\\j\not=k}}^n\Pi_k^{-1}Q_j^{-1}\gamma_k\beta_jb_k^{-3}b_j^{-3}K_b^{\left(2\right)}\left(\frac{X_i-X_k}{b_k}\right)K_b^{\left(2\right)}\left(\frac{X_i-X_j}{b_j}\right)Y_iY_j,\label{I2:rec}\\
\widehat{I}_3&=&\frac{\Pi_n^2}{n}\sum_{\substack{i,j,k,l=1\\j\not=k\not=l}}^n\Pi_j^{-1}\Pi_k^{-1}\gamma_j\gamma_kb_j^{-3}b_k^{-3}K_b^{\left(2\right)}\left(\frac{X_i-X_j}{b_j}\right)K_b^{\left(2\right)}\left(\frac{X_i-X_k}{b_k}\right)Y_iY_l,\label{I3:rec}\\
\widehat{I}_4&=&\frac{\Pi_n}{n}\sum_{\substack{i,k=1\\i\not=k}}^n\Pi_k^{-1}\gamma_kb_k^{-1}K_b\left(\frac{X_i-X_k}{b_k}\right)Y_i^2,\label{I4:rec}\\
\widehat{I}_5&=&\frac{Q_n}{n}\sum_{\substack{i,k=1\\i\not=k}}^nQ_k^{-1}\beta_kb_k^{-1}K_b\left(\frac{X_i-X_k}{b_k}\right)Y_iY_k,\label{I5:rec}
\end{eqnarray}
where $K_b$ is a kernel and $b_n$ is the associated bandwidth.\\
In practice, we take
\begin{eqnarray}\label{eq:h:initial}
b_n=n^{-\beta}\min\left\{\widehat{s},\frac{Q_3-Q_1}{1.349}\right\},\quad\beta \in \left]0,1\right[
\end{eqnarray} 
(see \citet{Sil86}) where $\widehat{s}$ the sample standard deviation, and $Q_1$, $Q_3$ denoting the first and third quartiles, respectively.\\
We followed the same steps as in~\citet{Sla14a} and we showed that in order to minimize the $MISE$ of $\widehat{I}_1$ respectively of $\widehat{I}_2$, $\widehat{I}_3$, $\widehat{I}_4$ and $\widehat{I}_5$, the pilot bandwidth $\left(b_n\right)$ must belong to $\mathcal{GS}\left(-3/14\right)$, respectively to $\mathcal{GS}\left(-3/14\right)$, $\mathcal{GS}\left(-3/14\right)$, $\mathcal{GS}\left(-2/5\right)$ and $\mathcal{GS}\left(-2/5\right)$.\\

Finally, the plug-in estimator of the bandwidth $\left(h_n\right)$ using the semi-recursive estimators defined in~(\ref{eq:rn}) with the stepsizes $\left(\gamma_n,\beta_n\right)=\left(n^{-1},n^{-1}\right)$. 
\begin{eqnarray}\label{hoptim:gamma:MISE}
\left(\left(\frac{3}{10}\right)^{1/5}\left(\frac{\widehat{I}_4-\widehat{I}_5}{\widehat{I}_1+\widehat{I}_3-2\widehat{I}_2}\right)^{1/5}\left\{\frac{R\left(K\right)}{\mu_2^2\left(K\right)}\right\}^{1/5}n^{-1/5}\right),
\end{eqnarray}
\begin{eqnarray*}
\widehat{MWISE}\left[r_n\right]&=&\frac{5}{4}\frac{1}{2^{4/5}}\left(\frac{5}{3}\right)^{6/5}\left(\widehat{I}_4-\widehat{I}_5\right)^{4/5}\times \left(\widehat{I}_1+\widehat{I}_3-2\widehat{I}_2\right)^{1/5}\Theta\left(K\right)n^{-4/5}\nonumber\\
&&+o\left(n^{-4/5}\right).
\end{eqnarray*}
Let us now consider the stepsize $\left(\beta_n\right)=\left(\left(1-a\right)n^{-1}\right)$, the case which minimize the variance of $a_n\left(x\right)$ combined with the stepsize $\left(\gamma_n\right)=\left(n^{-1}\right)$, the case which minimize the $MISE$ of $f_n$, it follows from~(\ref{MWISE:hoptim}), that
\begin{eqnarray}\label{MISE:gamma:var}
MWISE\left[r_n\right]&=&5^{1/5}\left(I_4-\frac{23}{24}I_5\right)^{4/5}\times \left(I_1+\frac{25}{36}I_3-\frac{5}{3}I_2\right)^{1/5}\Theta\left(K\right)n^{-4/5}\nonumber\\
&&+o\left(n^{-4/5}\right),
\end{eqnarray}
and from~(\ref{hoptim}), that the plug-in estimator of the bandwidth $\left(h_n\right)$ using the semi-recursive estimators defined in~(\ref{eq:rn}) is given by
\begin{eqnarray}\label{hoptim:gamma:var}
\left(\left(\frac{1}{5}\right)^{1/5}\left(\frac{\widehat{I}_4-\frac{23}{24}\widehat{I}_5}{\widehat{I}_1+\frac{25}{36}\widehat{I}_3-\frac{5}{3}\widehat{I}_2}\right)^{1/5}\left\{\frac{R\left(K\right)}{\mu_2^2\left(K\right)}\right\}^{1/5}n^{-1/5}\right),
\end{eqnarray}
and it follows from~(\ref{MWISE:hoptim}), that the plug-in $MWISE$ of the proposed estimator~(\ref{eq:rn}) using the stepsizes $\left(\gamma_n,\beta_n\right)=\left(n^{-1},\left(1-a\right)n^{-1}\right)$ is given by
\begin{eqnarray*}
\widehat{MWISE}\left[r_n\right]&=&5^{1/5}\left(\widehat{I}_4-\frac{23}{24}\widehat{I}_5\right)^{4/5}\times \left(\widehat{I}_1+\frac{25}{36}\widehat{I}_3-\frac{5}{3}\widehat{I}_2\right)^{1/5}\Theta\left(K\right)n^{-4/5}\nonumber\\
&&+o\left(n^{-4/5}\right).
\end{eqnarray*}
Let us now provide the case when the stepsize $\left(\gamma_n\right)$ is chosen to minimize the variance of $f_n$. It was shown in~\citet{Mok09a} and considered in~\citet{Sla13} that to minimize the variance of $f_n$, the stepsize $\left(\gamma_n\right)$ must be chosen in $\mathcal{GS}\left(-1\right)$ and must satisfy $\lim_{n\to\infty}n\gamma_n=1-a$. We consider here the case $\left(\gamma_n\right)=\left(\left(1-a\right)n^{-1}\right)$. Our first result is the following proposition, which gives the bias and the variance of $r_n$ in the special case of $\left(\gamma_n\right)=\left(\left(1-a\right)n^{-1}\right)$.
\begin{prop}[Bias and variance of $r_n$]\label{prop:bias:var:rn:bis}
Let Assumptions $\left(A1\right)-\left(A3\right)$ hold, and suppose that $\left(\gamma_n\right)=\left(\left(1-a\right)n^{-1}\right)$.
\begin{enumerate}
\item If $a\in ]0, \beta/5]$, then
\begin{eqnarray}\label{bias:rn:min:var}
\mathbb{E}\left[r_n\left(x\right)\right]-r\left(x\right)&=&\frac{1}{2f\left(x\right)}\left(\frac{a^{\left(2\right)}\left(x\right)}{\left(1-2a\xi\right)}-\frac{1-a}{\left(1-3a\right)}r\left(x\right)f^{\left(2\right)}\left(x\right)\right)h_n^2\mu_2\left(K\right)\nonumber\\
&&+o\left(h_n^2\right).
\end{eqnarray}
If $a\in  ]\beta/5, 1[$, then
\begin{eqnarray}\label{bias:rn:bis:min:var}
\mathbb{E}\left[r_n\left(x\right)\right]-r\left(x\right)=o\left(\sqrt{\beta_nh_n^{-1}}\right).
\end{eqnarray}
\item If $a\in [\beta/5, 1[$, then
\begin{eqnarray}
Var\left[r_n\left(x\right)\right]&=&\frac{\beta_n}{h_n}\left\{\frac{\mathbb{E}\left[Y^2\vert X=x\right]}{\left(2-\left(\beta-a\right)\xi\right)f\left(x\right)}-\left(1-a\right)\xi \frac{r^2\left(x\right)}{f\left(x\right)}\right\}R\left(K\right)\nonumber\\
&&+o\left(\frac{\beta_n}{h_n}\right).\label{var:rn:min:var}
\end{eqnarray}
If $a\in ]0,\beta/5[$, then
\begin{eqnarray}\label{var:rn:rep:min:var}
Var\left[r_n\left(x\right)\right]=o\left(h_n^4\right).
\end{eqnarray}
\item If $\lim_{n\to \infty}\left(n\beta_n\right)>\max\left\{2a, \left(a-\beta\right)/2\right\}$, then~(\ref{bias:rn:min:var}) and~(\ref{var:rn:min:var}) hold simultaneously.
\end{enumerate}
\end{prop}
The bias and the variance of the estimator $r_n$ defined by the stochastic approximation algorithm~(\ref{eq:rn}) then heavily depend on the choice of the stepsizes $\left(\gamma_n\right)$ and $\left(\beta_n\right)$.\\
Let us first state the following theorem, which gives the weak convergence rate of the estimator $r_n$ defined in~(\ref{eq:rn}).
\begin{theor}[Weak pointwise convergence rate]\label{theo:TLC2}
Let Assumptions $\left(A1\right)-\left(A3\right)$ hold, and suppose that $\left(\gamma_n\right)=\left(\left(1-a\right)n^{-1}\right)$.
\begin{enumerate}
\item If there exists $c\geq 0$ such that $\beta_n^{-1}h_n^5\to c$, then
\begin{eqnarray*}
\sqrt{\beta_n^{-1}h_n}\left(r_{n}\left(x\right)-r\left( x\right) \right) \stackrel{\mathcal{D}}{\rightarrow} &
\mathcal{N}\left(\sqrt{c}B_{a,\xi}^{\left(2\right)}\left(x\right),V_{a,\xi,\beta}^{\left(2\right)}\left(x\right)\right),
\end{eqnarray*}
where
\begin{eqnarray*}
B_{a,\xi}^{\left(2\right)}\left(x\right)&=&\frac{1}{2f\left(x\right)}\left(\frac{a^{\left(2\right)}\left(x\right)}{\left(1-2a\xi\right)}-\frac{\left(1-a\right)}{\left(1-3a\right)}r\left(x\right)f^{\left(2\right)}\left(x\right)\right)\mu_2\left(K\right),\\
V_{a,\xi,\beta}^{\left(2\right)}\left(x\right)&=&\left\{\frac{\mathbb{E}\left[Y^2\vert X=x\right]}{\left(2-\left(\beta-a\right)\xi\right)f\left(x\right)}-\left(1-a\right)\xi \frac{r^2\left(x\right)}{f\left(x\right)}\right\}R\left(K\right).
\end{eqnarray*}
\item If $nh_{n}^{5} \rightarrow \infty $, then  
\begin{eqnarray*}
\frac{1}{h_{n}^{2}}\left(r_{n}\left(x\right)-r\left(x\right) \right) \stackrel{\mathbb{P}}{\rightarrow } B_{a,\xi}^{\left(2\right)}\left(x\right).
\end{eqnarray*}
\end{enumerate}
\end{theor}
The following corollary gives the weak convergence rate of $r_n$ in the two special cases;  $\left(\gamma_n,\beta_n\right)=\left(\left(1-a\right)n^{-1},n^{-1}\right)$ and $\left(\gamma_n,\beta_n\right)=\left(\left(1-a\right)n^{-1},\left(1-a\right)n^{-1}\right)$ respectively. 
\begin{coro}[Weak pointwise convergence rate]\label{coro:TLC2}
Let Assumptions $\left(A1\right)-\left(A3\right)$ hold. 
\begin{enumerate}
\item If we suppose that the stepsizes $\left(\gamma_n,\beta_n\right)=\left(\left(1-a\right)n^{-1},n^{-1}\right)$, and if there exists $c\geq 0$ such that $nh_n^5\to c$, then
\begin{eqnarray*}
\sqrt{nh_n}\left(r_{n}\left(x\right)-r\left( x\right) \right) \stackrel{\mathcal{D}}{\rightarrow} &
\mathcal{N}\left(\sqrt{c}B_{a,1}^{\left(2\right)}\left(x\right),V_{a,1,1}^{\left(2\right)}\left(x\right)\right).
\end{eqnarray*}
\item If we suppose that the stepsizes $\left(\gamma_n,\beta_n\right)=\left(\left(1-a\right)n^{-1},\left(1-a\right)n^{-1}\right)$, and if there exists $c\geq 0$ such that $nh_n^5\to c$, then
\begin{eqnarray*}
\sqrt{nh_n}\left(r_{n}\left(x\right)-r\left( x\right) \right) \stackrel{\mathcal{D}}{\rightarrow} &
\mathcal{N}\left(\sqrt{c}B_{a,\left(1-a\right)^{-1}}^{\left(2\right)}\left(x\right),V_{a,\left(1-a\right)^{-1},1}^{\left(2\right)}\left(x\right)\right).
\end{eqnarray*}
\end{enumerate}
\end{coro} 

The following proposition gives the $MWISE$ of $r_n$ in the case when $\left(\gamma_n\right)$ is chosen to minimize the variance of $f_n$.
\begin{prop}[$MWISE$ of $r_n$]\label{prop:MISE:var:rn}
Let Assumptions $\left(A1\right)-\left(A3\right)$ hold, and suppose that $\left(\gamma_n\right)=\left(\left(1-a\right)n^{-1}\right)$.
\begin{enumerate}
\item If $a\in  ]0, \beta/5[$, then
\begin{eqnarray*}
MWISE\left[r_n\right]&=&\frac{1}{4}\left(\frac{I_1}{\left(1-2a\xi\right)^2}+\frac{\left(1-a\right)^2}{\left(1-3a\right)^2}I_3-2\frac{\left(1-a\right)}{\left(1-3a\right)\left(1-2a\xi\right)}I_2\right)h_n^4\mu_2^2\left(K\right)\\
&&+o\left(h_n^4\right).
\end{eqnarray*}
\item If $a=\beta/5$, then
\begin{eqnarray*}
MWISE\left[r_n\right]&=&\frac{\beta_n}{h_n}\left(\frac{I_4}{\left(2-\left(\beta-a\right)\xi\right)}-\left(1-a\right)\xi I_5\right)R\left(K\right)\\
&&+\frac{1}{4}\left(\frac{I_1}{\left(1-2a\xi\right)^2}+\frac{\left(1-a\right)^2}{\left(1-3a\right)^2}I_3-2\frac{\left(1-a\right)}{\left(1-3a\right)\left(1-2a\xi\right)}I_2\right)h_n^4\mu_2^2\left(K\right)\\
&&+o\left(h_n^4\right).
\end{eqnarray*}
\item If $a\in ]\beta/5, 1[$, then
\begin{eqnarray*}
MWISE\left[r_n\right]&=&\frac{\beta_n}{h_n}\left(\frac{I_4}{\left(2-\left(\beta-a\right)\xi\right)}-\left(1-a\right)\xi I_5\right)R\left(K\right)+o\left(\frac{\beta_n}{h_n}\right).
\end{eqnarray*}
\end{enumerate}
\end{prop}
The following corollary ensures that the bandwidth which minimize the $MWISE$ depend on the stepsize $\left(\beta_n\right)$ and then the corresponding $MWISE$ depend also on the stepsize $\left(\beta_n\right)$.
\begin{coro}\label{coro:hn:MISE:var}
Let Assumptions $\left(A1\right)-\left(A3\right)$ hold, and suppose that $\left(\gamma_n\right)=\left(\left(1-a\right)n^{-1}\right)$. To minimize the $MWISE$ of $r_n$, the stepsize $\left(\beta_n\right)$ must be chosen in $\mathcal{GS}\left(-1\right)$, the bandwidth $\left(h_n\right)$ must equal 
\begin{eqnarray*}
\left(\left\{\frac{\frac{I_4}{\left(2-\left(\beta-a\right)\xi\right)}-\left(1-a\right)\xi I_5}{\frac{I_1}{\left(1-2a\xi\right)^2}+\frac{\left(1-a\right)^2}{\left(1-3a\right)^2}I_3-2\frac{\left(1-a\right)}{\left(1-3a\right)\left(1-2a\xi\right)}I_2}
\right\}^{1/5}\left\{\frac{R\left(K\right)}{\mu_2^2\left(K\right)}\right\}^{1/5}\beta_n^{1/5}\right).
\end{eqnarray*}
Then, we have
\begin{eqnarray*} MWISE\left[r_n\right]&=&\frac{5}{4}\left(\frac{I_4}{\left(2-\left(\beta-a\right)\xi\right)}-\left(1-a\right)\xi I_5\right)^{4/5}\\
&&\times\left(\frac{I_1}{\left(1-2a\xi\right)^2}+\frac{\left(1-a\right)^2}{\left(1-3a\right)^2}I_3-2\frac{\left(1-a\right)}{\left(1-3a\right)\left(1-2a\xi\right)}I_2\right)^{1/5}\Theta\left(K\right)\beta_n^{4/5}\\
&&+o\left(\beta_n^{4/5}\right).
\end{eqnarray*}
\end{coro}
The following corollary shows that, for a special choice of the stepsize $\left(\beta_n\right)=\left(\beta_0n^{-1}\right)$, which fulfilled that $\lim_{n\to \infty}n\beta_n=\beta_0$ and that $\left(\beta_n\right)\in \mathcal{GS}\left(-1\right)$, the optimal value for $h_n$ depend on $\beta_0$ and then the corresponding $MWISE$ depend on $\beta_0$.
\begin{coro}\label{coro:hn:MISE:var:bis}
Let Assumptions $\left(A1\right)-\left(A3\right)$ hold, and suppose that $\left(\gamma_n\right)=\left(\left(1-a\right)n^{-1}\right)$. To minimize the $MWISE$ of $r_n$, the stepsize $\left(\beta_n\right)$ must be chosen in $\mathcal{GS}\left(-1\right)$, $\lim_{n\to \infty}n\beta_n=\beta_0$, the bandwidth $\left(h_n\right)$ must equal 
\begin{eqnarray}\label{hoptim:beta}
\left(\left(\frac{\beta_0-2/5}{2}\right)^{1/5}
\left(\frac{I_4-\frac{8}{5}\frac{\left(\beta_0-2/5\right)}{\beta_0^2}I_5}{I_1+4\left(\frac{\beta_0-2/5}{\beta_0}\right)^2I_3-4\left(\frac{\beta_0-2/5}{\beta_0}\right)I_2}\right)^{1/5}
\left\{\frac{R\left(K\right)}{\mu_2^2\left(K\right)}\right\}^{1/5}n^{-1/5}\right),
\end{eqnarray}
and we then have
\begin{eqnarray}\label{MWISE:hoptim:beta} MWISE\left[r_n\right]&=&\frac{5}{4}\frac{1}{2^{4/5}}\frac{\beta_0^2}{\left(\beta_0-2/5\right)^{6/5}}\left(I_4-\frac{8}{5}\frac{\beta_0-2/5}{\beta_0^2}I_5\right)^{4/5}\nonumber\\
&&\times \left(I_1+4\left(\frac{\beta_0-2/5}{\beta_0}\right)^2I_3-4\left(\frac{\beta_0-2/5}{\beta_0}\right)I_2\right)^{1/5} \Theta\left(K\right)n^{-4/5}\nonumber\\
&&+o\left(n^{-4/5}\right).
\end{eqnarray}
Moreover, the minimum of $\beta_0^2\left(\beta_0-2/5\right)^{-6/5}$ is reached at $\beta_0=1$, then the bandwidth $\left(h_n\right)$ must equal 
\begin{eqnarray*}
\left(\left(\frac{3}{10}\right)^{1/5}\left(\frac{I_4-\frac{24}{25}I_5}{I_1+\frac{36}{25}I_3-\frac{12}{5}I_2}\right)^{1/5}\left\{\frac{R\left(K\right)}{\mu_2^2\left(K\right)}\right\}^{1/5}n^{-1/5}\right),
\end{eqnarray*}
and we then have
\begin{eqnarray}\label{MISE:beta:MISE} MWISE\left[r_n\right]&=&\frac{5}{4}\frac{1}{2^{4/5}}\left(\frac{5}{3}\right)^{6/5}\left(I_4-\frac{24}{25}I_5\right)^{4/5}\times \left(I_1+\frac{36}{25}I_3-\frac{12}{5}I_2\right)^{1/5}\Theta\left(K\right)n^{-4/5}\nonumber\\
&&+o\left(n^{-4/5}\right).
\end{eqnarray}
\end{coro}
In order to estimate the optimal bandwidth (\ref{hoptim:gamma}), we must estimate $I_1$, $I_2$, $I_3$, $I_4$ and $I_5$. We use the kernel estimators defined in~(\ref{I1:rec}),~(\ref{I2:rec}),~(\ref{I3:rec}),~(\ref{I4:rec}) and (\ref{I5:rec}).
We showed that in order to minimize the $MISE$ of $\widehat{I}_1$ respectively of $\widehat{I}_2$, $\widehat{I}_3$, $\widehat{I}_4$ and $\widehat{I}_5$, the pilot bandwidth $\left(b_n\right)$ must belong to $\mathcal{GS}\left(-3/14\right)$, respectively to $\mathcal{GS}\left(-3/14\right)$, $\mathcal{GS}\left(-3/14\right)$, $\mathcal{GS}\left(-2/5\right)$ and $\mathcal{GS}\left(-2/5\right)$.\\
Finally, the plug-in estimator of the bandwidth $\left(h_n\right)$ using the semi-recursive estimators defined in~(\ref{eq:rn}) with the stepsizes $\left(\gamma_n,\beta_n\right)=\left(\left(1-a\right)n^{-1},n^{-1}\right)$. 
\begin{eqnarray}\label{hoptim:beta:MISE}
\left(\left(\frac{3}{10}\right)^{1/5}\left(\frac{\widehat{I}_4-\frac{24}{25}\widehat{I}_5}{\widehat{I}_1+\frac{36}{25}\widehat{I}_3-\frac{12}{5}\widehat{I}_2}\right)^{1/5}\left\{\frac{R\left(K\right)}{\mu_2^2\left(K\right)}\right\}^{1/5}n^{-1/5}\right),
\end{eqnarray}
\begin{eqnarray*}
\widehat{MWISE}\left[r_n\right]&=&\frac{5}{4}\frac{1}{2^{4/5}}\left(\frac{5}{3}\right)^{6/5}\left(\widehat{I}_4-\frac{24}{25}\widehat{I}_5\right)^{4/5}\times \left(\widehat{I}_1+\frac{36}{25}\widehat{I}_3-\frac{12}{5}\widehat{I}_2\right)^{1/5}\Theta\left(K\right)n^{-4/5}\nonumber\\
&&+o\left(n^{-4/5}\right).
\end{eqnarray*}
Let us now consider the stepsize $\left(\beta_n\right)=\left(\left(1-a\right)n^{-1}\right)$, the case which minimize the variance of $a_n\left(x\right)$ combined with the stepsize $\left(\gamma_n\right)=\left(\left(1-a\right)n^{-1}\right)$, the case which minimize the variance of $f_n$, it follows from~(\ref{MWISE:hoptim:beta}), that
\begin{eqnarray}\label{MISE:beta:var}
MWISE\left[r_n\right]&=&5^{1/5}\left(I_4-I_5\right)^{4/5}\times \left(I_1+I_3-2I_2\right)^{1/5}\Theta\left(K\right)n^{-4/5}+o\left(n^{-4/5}\right),
\end{eqnarray}
and from (\ref{hoptim:beta}), that the plug-in estimator of the bandwidth $\left(h_n\right)$ using the semi-recursive estimators defined in~(\ref{eq:rn}) is given by 
\begin{eqnarray}\label{hoptim:beta:var}
\left(\left(\frac{1}{5}\right)^{1/5}\left(\frac{\widehat{I}_4-\widehat{I}_5}{\widehat{I}_1+\widehat{I}_3-2\widehat{I}_2}\right)^{1/5}\left\{\frac{R\left(K\right)}{\mu_2^2\left(K\right)}\right\}^{1/5}n^{-1/5}\right),
\end{eqnarray}
and it follows from~(\ref{MWISE:hoptim:beta}), that the plug-in $MWISE$ of the proposed estimator~(\ref{eq:rn}) using the stepsizes $\left(\gamma_n,\beta_n\right)=\left(\left(1-a\right)n^{-1},\left(1-a\right)n^{-1}\right)$, is given by 
\begin{eqnarray*}
\widehat{MWISE}\left[r_n\right]&=&5^{1/5}\left(\widehat{I}_4-\widehat{I}_5\right)^{4/5}\times \left(\widehat{I}_1+\widehat{I}_3-2\widehat{I}_2\right)^{1/5}\Theta\left(K\right)n^{-4/5}+o\left(n^{-4/5}\right).
\end{eqnarray*}

Now, let us recall that the bias and variance of Nadaraya-Watson's estimator $\widetilde{r}_n$ are given by
\begin{eqnarray*}
\mathbb{E}\left[\widetilde{r}_n\left(x\right)\right]-r\left(x\right)=\frac{1}{2}\left(a^{\left(2\right)}\left(x\right)-r\left(x\right)f^{\left(2\right)}\left(x\right)\right)f^{-1}\left(x\right)h_n^2\mu_2\left(K\right)+o\left(h_n^2\right),
\end{eqnarray*}
and
\begin{eqnarray*}
Var\left[\widetilde{r}_n\left(x\right)\right]&=&\frac{1}{nh_n}Var\left[Y\vert X=x\right]f^{-1}\left(x\right)R\left(K\right)+o\left(\frac{1}{nh_n}\right).
\end{eqnarray*}
It follows that,
\begin{eqnarray*} MWISE\left[\widetilde{r}_n\right]&=&\frac{1}{nh_n}\left(I_4-I_5\right)R\left(K\right)+\frac{1}{4}\left(I_1+I_3-2I_2\right)h_n^4\mu_2^2\left(K\right)+o\left(h_n^4+\frac{1}{nh_n}\right).
\end{eqnarray*}
Then, to minimize the $MWISE$ of $\widetilde{r}_n$, the bandwidth $\left(h_n\right)$ must equal to
\begin{eqnarray}\label{hoptim:rose}
\left(\left(\frac{I_4-I_5}{I_1+I_3-2I_2}\right)^{1/5}\left\{\frac{R\left(K\right)}{\mu_2^2\left(K\right)}\right\}^{1/5}n^{-1/5}\right),
\end{eqnarray}
and we have
\begin{eqnarray}\label{MISE:rose} MWISE\left[\widetilde{r}_n\right]=\frac{5}{4}\left(I_4-I_5\right)^{4/5}\left(I_1+I_3-2I_2\right)^{1/5}\Theta\left(K\right)n^{-4/5}
+o\left(n^{-4/5}\right).
\end{eqnarray}
To estimate the optimal bandwidth~(\ref{hoptim:rose}), we must estimate $I_1$, $I_2$, $I_3$, $I_4$ and $I_5$. We use the following kernel estimator of $I_1$, $I_2$, $I_3$, $I_4$ and $I_5$:
\begin{eqnarray*}
\widetilde{I}_1&=&\frac{1}{n^3b_n^{6}}\sum_{\substack{i,j,k=1\\j\not=k}}^nK_b^{\left(2\right)}\left(\frac{X_i-X_j}{b_n}\right)K_b^{\left(2\right)}\left(\frac{X_i-X_k}{b_n}\right)Y_jY_k,\\
\widetilde{I}_2&=&\frac{1}{n^3b_n^{6}}\sum_{\substack{i,j,k=1\\j\not=k}}^nK_b^{\left(2\right)}\left(\frac{X_i-X_j}{b_n}\right)K_b^{\left(2\right)}\left(\frac{X_i-X_k}{b_n}\right)Y_iY_j,\\
\widetilde{I}_3&=&\frac{1}{n^4b_n^{6}}\sum_{\substack{i,j,k,l=1\\j\not=k\not=l}}^nK_b^{\left(2\right)}\left(\frac{X_i-X_j}{b_n}\right)K_b^{\left(2\right)}\left(\frac{X_i-X_k}{b_n}\right)Y_iY_l,\\
\widetilde{I}_4&=&\frac{1}{n^2b_n}\sum_{\substack{i,j=1\\i\not=j}}^nK_b\left(\frac{X_i-X_j}{b_n}\right)Y_i^2,\\
\widetilde{I}_5&=&\frac{1}{n^2b_n}\sum_{\substack{i,k=1\\i\not=k}}^nK_b\left(\frac{X_i-X_k}{b_n}\right)Y_iY_k,
\end{eqnarray*}
where $K_b$ is a kernel and $b_n$ is the associated bandwidth given in~(\ref{eq:h:initial}).\\ 
We showed that in order to minimize the $MISE$ of $\widetilde{I}_1$ respectively of $\widetilde{I}_2$, $\widetilde{I}_3$, $\widetilde{I}_4$ and $\widetilde{I}_5$, the pilot bandwidth $\left(b_n\right)$ must belong to $\mathcal{GS}\left(-3/14\right)$, respectively to $\mathcal{GS}\left(-3/14\right)$, $\mathcal{GS}\left(-3/14\right)$, $\mathcal{GS}\left(-2/5\right)$ and $\mathcal{GS}\left(-2/5\right)$.\\
Then the plug-in estimator of the bandwidth $\left(h_n\right)$ using the nonrecursive estimator~(\ref{eq:Nad}), is given by
\begin{eqnarray}\label{hoptim:rose:plug:in}
\left(\left(\frac{\widetilde{I}_4-\widetilde{I}_5}{\widetilde{I}_1+\widetilde{I}_3-2\widetilde{I}_2}\right)^{1/5}\left\{\frac{R\left(K\right)}{\mu_2^2\left(K\right)}\right\}^{1/5}n^{-1/5}\right),
\end{eqnarray}
and the plug-in of the $MWISE$ of the nonrecursive estimator~(\ref{eq:Nad}), is given by
\begin{eqnarray*}
\widetilde{MWISE}\left[\widetilde{r}_n\right]=\frac{5}{4}\left(\widetilde{I}_4-\widetilde{I}_5\right)^{4/5}\left(\widetilde{I}_1+\widetilde{I}_3-2\widetilde{I}_2\right)^{1/5}\Theta\left(K\right)n^{-4/5}
+o\left(n^{-4/5}\right).
\end{eqnarray*}
Finally, it follows from~(\ref{MISE:gamma:MISE}),~(\ref{MISE:gamma:var}),~(\ref{MISE:beta:MISE}),~(\ref{MISE:beta:var}) and~(\ref{MISE:rose}), that:
\begin{description}
\item The $MWISE$ of the proposed estimator~(\ref{eq:rn}) with the choice of the stepsizes $\left(\gamma_n,\beta_n\right)=\left(n^{-1},n^{-1}\right)$ is $1.06$ larger than the nonrecursive estimator~(\ref{eq:Nad}).
\item The $MWISE$ of the proposed estimator~(\ref{eq:rn}) with the choice of the stepsizes $\left(\gamma_n,\beta_n\right)=\left(\left(1-a\right)n^{-1},\left(1-a\right)n^{-1}\right)$ is $1.1$ larger than the nonrecursive estimator~(\ref{eq:Nad}).
\item We can't compare the $MWISE$ of the proposed estimator~(\ref{eq:rn}) with the choice of the stepsizes $\left(\gamma_n,\beta_n\right)=\left(n^{-1},\left(1-a\right)n^{-1}\right)$ (respectively, the $MWISE$ of the proposed estimator~(\ref{eq:rn}) with the choice of the stepsizes $\left(\gamma_n,\beta_n\right)=\left(\left(1-a\right)n^{-1},n^{-1}\right)$) neither to the $MWISE$ of the others proposed estimators nor to the $MWISE$ of the nonrecursive estimator~(\ref{eq:Nad}). 
\end{description}
\section{Applications}\label{section:app}
The aim of our applications is to compare the performance of the semi-recursive estimators defined in~(\ref{eq:rn}) with that of the nonrecursive Nadaraya-Watson's estimator defined in~(\ref{eq:Nad}).
\begin{description}
\item When applying $r_n$ one need to choose three quantities:
\begin{itemize}
\item The function $K$, we choose the Normal kernel. 
\item The stepsizes $\left(\gamma_n,\beta_n\right)$ equal respectively to $\left(n^{-1},n^{-1}\right)$, $\left(n^{-1},\left(1-a\right)n^{-1}\right)$, \\
$\left(\left(1-a\right)n^{-1},n^{-1}\right)$ or $\left(\left(1-a\right)n^{-1},\left(1-a\right)n^{-1}\right)$. These four choices are referred to as \texttt{Recursive 1}, \texttt{2}, \texttt{3} and \texttt{4} respectively.
\item The bandwidth $\left(h_n\right)$ is chosen to be equal respectively to~(\ref{hoptim:gamma:MISE}) for (\texttt{Recursive 1}), (\ref{hoptim:gamma:var}) for (\texttt{Recursive 2}), (\ref{hoptim:beta:var}) for (\texttt{Recursive 3}) and (\ref{hoptim:beta:MISE}) for (\texttt{Recursive 4}).
\end{itemize}
\item When applying $\widetilde{r}_n$ one need to choose two quantities: 
\begin{itemize}
\item The function $K$, as in the semi-recursive framework, we use the Normal kernel. 
\item The bandwidth $\left(h_n\right)$ is chosen to be equal to~(\ref{hoptim:rose:plug:in}). 
\end{itemize}
\end{description}

\subsection{Simulations}\label{subsection:simu}
Througthout this subsection, we consider the regression model
\begin{eqnarray*}
Y=r\left(X\right)+\varepsilon,
\end{eqnarray*}
where $X$ is $\mathcal{N}\left(0,1\right)$-distributed and $\varepsilon$ is $\mathcal{N}\left(0,\sigma\right)$-distributed, with $\sigma$ is chosen in the interval $\left[0.1,2\right]$.\\
In order to investigate the comparison between the proposed estimators, we consider two regression functions : cosine function $r\left(x\right)=\cos\left(x\right)$ (see Table~\ref{Tab:1}) and the following function $r\left(x\right)=\left(1+\exp\left(x\right)\right)^{-1}$ (see Table~\ref{Tab:2}). For each fixed $\sigma\in \left[0.1,2\right]$, the number of simulations is $500$. We denote by $r_i^*$ the true regression function, and by $r_i$ the considered regression estimators, and then we compute the Mean Squared Error ($MSE=n^{-1}\sum_i\left(r_i-r_i^*\right)^2$). 

\paragraph{Computational cost}
The advantage of recursive estimators on their nonrecursive version is that their update, from a sample of size $n$ to one of size $n+1$, require less computations. 
Performing all the proposed methods, we report the total \texttt{CPU} time values for each considered regression function and for each fixed $\sigma$ and for each sample size in Tables~\ref{Tab:1} and~\ref{Tab:2}. For the two tables we give the \texttt{CPU} time in seconds.
\begin{table}[!h]
\begin{eqnarray*}
\begin{tabular}{lccccccc}
& Nadaraya  & Recursive $1$  & Recursive $2$  & Recursive $3$ & Recursive $4$\\ \hline
$n=100$&&  &$\sigma=0.1$  &&& \\
$MSE$ &  $0.000812$ & $ 0.000748$ & $0.000764$ & ${\bf 0.000567}$ & $0.000667$\\
$\texttt{CPU}$ & $238$ & $184$ & $170$& ${\bf 154}$ & $164$\\
$n=200$&&  & &&& \\
$MSE$ &  $0.000507$ & $0.000483$ & $0.000508$ & ${\bf 0.000366}$ & $0.000443$\\
$\texttt{CPU}$ & $835$ & $514$ & $509$ & ${\bf 464}$& $470$\\
$n=500$&&  &  &&&\\
$MSE$ &  $0.000284$ & $0.000279$ & $0.000294$ & ${\bf 0.000217}$ & $0.000260$\\
$\texttt{CPU}$ & $3679$ & $2185$ & $1973$ & $1966$ & ${\bf 1865}$\\
\hline
$n=100$&&  &$\sigma=0.5$  &&& \\
$MSE$ &  $0.004486$ & $0.004447$ & $0.004286$ & ${\bf 0.003729}$ & $0.004184$\\
$\texttt{CPU}$ & $231$ & $143$ & $135$ & $137$ & ${\bf 129}$\\
$n=200$&&  &  &&& \\
$MSE$ &  $0.002331$ & $0.002337$ & $0.002142$ & ${\bf 0.001929}$ & $0.002141$\\
$\texttt{CPU}$ & $885$ & $568$ & $549$ & $485$& ${\bf 457}$\\
$n=500$&&  &  &&&\\
$MSE$ &  $0.001372$ & $0.001411$ & $0.001265$ & ${\bf 0.001174}$ & $0.001291$\\
$\texttt{CPU}$ & $3498$ & $2049$ & ${\bf 1943}$ & $2242$ & $2045$\\
\hline
$n=100$&&  &$\sigma=1$  &&& \\
$MSE$ &  ${\bf 0.013960}$ & $0.021204$ & $0.020982$ & $0.021476$ & $0.021832$\\
$\texttt{CPU}$ & $246$ & $166$ & ${\bf 136}$& $146$ & $137$\\
$n=200$&&  &  &&& \\
$MSE$ &  ${\bf 0.006016}$ & $0.010935$ & $0.008714$ & $0.012524$ & $0.011657$\\
$\texttt{CPU}$ & $831$ & $580$ & $519$ & $541$ & ${\bf 505}$\\
$n=500$&&  &  &&&\\
$MSE$ &  $0.001916$ & ${\bf 0.001816}$ & $0.002268$ & $0.003018$ & $0.001972$\\
$\texttt{CPU}$ & $3801$ & $2193$ &$2043$ & $2024$ & ${\bf 1875}$\\
\hline
\end{tabular}
\end{eqnarray*}
\caption{Quantitative comparison between the nonrecursive estimator~(\ref{eq:Nad}) and four recursive estimators; recursive $1$ correspond to the estimator~(\ref{eq:rn}) with the choice $\left(\gamma_n,\beta_n\right)=\left(n^{-1},n^{-1}\right)$, recursive $2$ correspond to the estimator~(\ref{eq:rn}) with the choice $\left(\gamma_n,\beta_n\right)=\left(n^{-1},\left(1-a\right)n^{-1}\right)$, recursive $3$ correspond to the estimator~(\ref{eq:rn}) with the choice $\left(\gamma_n,\beta_n\right)=\left(\left(1-a\right)n^{-1},n^{-1}\right)$ and recursive $4$ correspond to the estimator~(\ref{eq:rn}) with the choice $\left(\gamma_n,\beta_n\right)=\left(\left(1-a\right)n^{-1},\left(1-a\right)n^{-1}\right)$. Here we consider the regression function $r\left(x\right)=\cos\left(x\right)$, $X\sim \mathcal{N}\left(0,1\right)$ and $\varepsilon\sim \mathcal{N}\left(0,\sigma\right)$ with $\sigma=0.1$ in the first block, $\sigma=0.5$ in the second block and $\sigma=1$ in the last block, we consider three sample sizes $n=100$, $n=200$ and $n=500$, the number of simulations is $500$, and we compute the Mean squared error ($MSE$) and the \texttt{CPU} time in seconds.}\label{Tab:1}
\end{table}

\begin{table}[!h]
\begin{eqnarray*}
\begin{tabular}{lccccccc}
& Nadaraya  & Recursive $1$  & Recursive $2$  & Recursive $3$ & Recursive $4$\\ \hline
$n=100$&&  &$\sigma=0.1$  &&& \\
$MSE$ &  $1.31e^{-04}$ & $1.15e^{-04}$ & ${\bf 6.22e^{-05}}$ & $1.71e^{-04}$ & $ 1.03e^{-04}$\\
$\texttt{CPU}$ & $249$ & $184$ & ${\bf 135}$& $146$ & $146$\\
$n=200$&&  & &&& \\
$MSE$ &  $4.38e^{-05}$ & $3.87e^{-05}$ & ${\bf 1.50e^{-05}}$ & $8.03e^{-05}$ & $3.63e^{-05}$\\
$\texttt{CPU}$ & $909$ & $524$ & $475$& $601$ & ${\bf 458}$\\
$n=500$&&  &  &&&\\
$MSE$ &  $5.70e^{-06}$ & $5.02e^{-06}$ & ${\bf 3.20e^{-06}}$ & $2.32e^{-05}$ & $4.29e^{-06}$\\
$\texttt{CPU}$ & $3708$ & $1803$ & $1672$& $1855$ & ${\bf 1483}$\\
\hline
$n=100$&&  &$\sigma=0.5$  &&& \\
$MSE$ &  $0.000351$ & $0.000325$ & ${\bf 0.000252}$ & $0.000350$ & $0.000296$\\
$\texttt{CPU}$ & $256$ & $144$ & $132$& $134$ & ${\bf 125}$\\
$n=200$&&  &  &&& \\
$MSE$ &  $0.000189$ & $0.000171$ & $0.000154$ & $0.000163$ & ${\bf 0.000151}$\\
$\texttt{CPU}$ & $873$ & $524$ & $483$& $576$ & ${\bf 451}$\\
$n=500$&&  &  &&&\\
$MSE$ &  $2.30e^{-05}$ & $2.25e^{-05}$ & $2.42e^{-05}$ & $3.351e^{-05}$ & ${\bf 2.06e^{-05}}$\\
$\texttt{CPU}$ & $4389$ & $2113$ & $1987$& $1999$ & ${\bf 1973}$\\
\hline
$n=100$&&  &$\sigma=2$  &&& \\
$MSE$ &  $0.003447$ & $0.003294$ & $0.003155$ & ${\bf 0.003132}$ & $0.003137$\\
$\texttt{CPU}$ & $294$ & $155$ & $173$& ${\bf 143}$ & $148$\\
$n=200$&&  &  &&& \\
$MSE$ &  $0.000160$ & $0.000152$ & $0.000162$ & ${\bf 0.000111}$ & $0.000189$\\
$\texttt{CPU}$ & $917$ & $503$ & $581$& $515$ & ${\bf 477}$\\
$n=500$&&  &  &&&\\
$MSE$ &  $6.56e^{-05}$ & $7.03e^{-05}$ & ${\bf 5.01e^{-05}}$ & $6.70e^{-05}$ & $5.39e^{-05}$\\
$\texttt{CPU}$ & $3643$ & $2105$ & $1951$& $1947$ & ${\bf 1877}$\\
\hline
\end{tabular}
\end{eqnarray*}
\caption{Quantitative comparison between the nonrecursive estimator~(\ref{eq:Nad}) and four recursive estimators; recursive $1$ correspond to the estimator~(\ref{eq:rn}) with the choice $\left(\gamma_n,\beta_n\right)=\left(n^{-1},n^{-1}\right)$, recursive $2$ correspond to the estimator~(\ref{eq:rn}) with the choice $\left(\gamma_n,\beta_n\right)=\left(n^{-1},\left(1-a\right)n^{-1}\right)$, recursive $3$ correspond to the estimator~(\ref{eq:rn}) with the choice $\left(\gamma_n,\beta_n\right)=\left(\left(1-a\right)n^{-1},n^{-1}\right)$ and recursive $4$ correspond to the estimator~(\ref{eq:rn}) with the choice $\left(\gamma_n,\beta_n\right)=\left(\left(1-a\right)n^{-1},\left(1-a\right)n^{-1}\right)$. Here we consider the regression function $r\left(x\right)=\left(1+\exp\left(x\right)\right)^{-1}$, $X\sim \mathcal{N}\left(0,1\right)$ and $\varepsilon\sim \mathcal{N}\left(0,\sigma\right)$ with $\sigma=0.1$ in the first block, $\sigma=0.5$ in the second block and $\sigma=2$ in the last block, we consider three sample sizes $n=100$, $n=200$ and $n=500$, the number of simulations is $500$, and we compute the Mean squared error ($MSE$) and the \texttt{CPU} time in seconds.}\label{Tab:2}
\end{table}

\subsection{Real Dataset}\label{subsection:real}
The CO2 dataset was available in the R package \texttt{Stat2Data} and contained $237$ observations on the following two variables; Day and CO2, for more details see~the station information system (GAWSIS). Scientists at a research station in Brotjacklriegel, Germany recorded CO2 levels, in parts per million, in the atmosphere for each day from the start of April through November in 2001.
\begin{center}
\begin{figure}[!h]
\includegraphics[width=0.6\textwidth,angle=270,clip=true,trim=40 0 0 0]{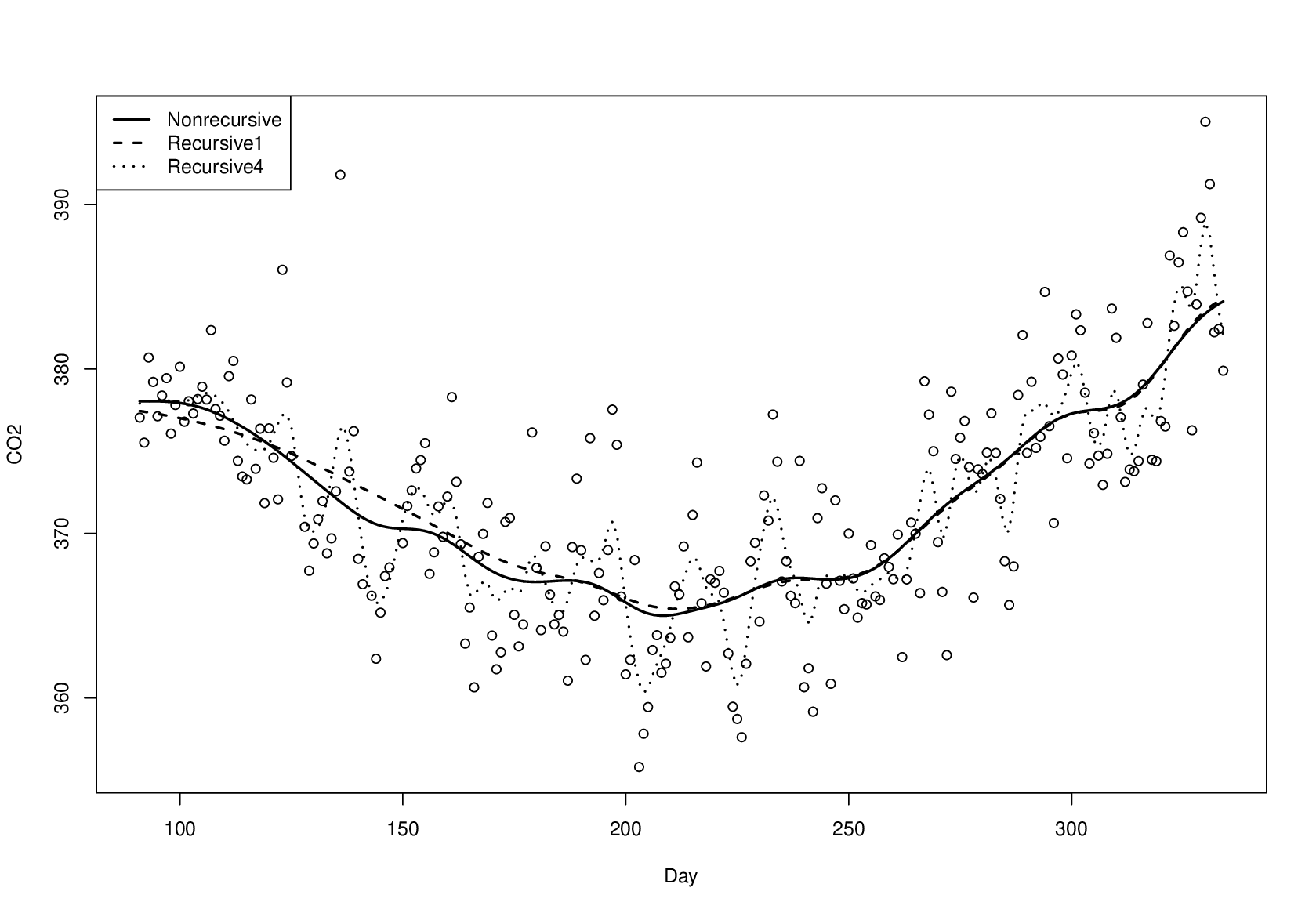}
\caption{The daily carbon dioxide measurements data with automatically bandwidth selection using the nonrecursive Nadaraya's estimator~(\ref{eq:Nad}) and two semi-recursive estimators~(\ref{eq:rn}) (\texttt{Recursive 1} and \texttt{Recursive 4}).}
\label{Fig:1}
\end{figure}
\end{center}
Figure~\ref{Fig:1} and Tables~\ref{Tab:1} and~\ref{Tab:2} indicate that
\begin{itemize}
\item The \texttt{Recursive 1} is very close to the nonrecursive estimator~(\ref{eq:Nad}).
\item The two estimators \texttt{Recursive 2} and \texttt{Recursive 3} can be better than the others estimators in many situations.
\item The \texttt{CPU} time are always faster using the proposed semi-recursive estimators and the reduction of \texttt{CPU} time goes from a minimum of $22.3\%$ to a maximum of $60\%$ compared to the nonrecursive estimator.
\end{itemize}

\section{Conclusion}\label{section:conclusion}

This paper propose an automatic selection of the bandwidth of the semi-recursive kernel estimators of a regression function defined by the stochastic approximation algorithm~(\ref{eq:rn}). The proposed estimators asymptotically follows normal distribution. The estimators are compared to the nonrecursive Nadaraya-Watson's regression estimator. We showed that, using some selected bandwidth and some particularly stepsizes, the proposed semi-recursive estimators will be very competitive to the nonrecursive one. The simulation study confirms the nice features of our proposed semi-recursive estimators and statisfactory improvement in the \texttt{CPU} time in comparison to the nonrecursive estimator. 

In conclusion, the proposed method allowed us to obtain quite similar results as the nonrecursive estimator proposed by~\citet{Nad64} and~\citet{Wat64}. Moreover, we plan to make an extensions of our method in future and to consider the case of the averaged R\'ev\'esz's regression estimators (see~\citet{Mok09b} and~\citet{Sla15a,Sla15b}) and the case of time series as in~\citet{Har90} in recursive way (see~\citet{Hua14}). 
\appendix
\section{Proofs} \label{section:proofs}
Throughout this section we use the following notation:
\begin{eqnarray*}
Q_n=\prod_{j=1}^{n}\left(1-\beta_j\right),\quad \Pi_n=\prod_{j=1}^{n}\left(1-\gamma_j\right),\quad \zeta_n=\Pi_nQ_n^{-1},
\end{eqnarray*}
\begin{eqnarray}
&&W_n\left(x\right)=h_n^{-1}K\left(\frac{x-X_n}{h_n}\right).\label{eq:12}\\
&&Z_n\left(x\right)=h_n^{-1}Y_nK\left(\frac{x-X_n}{h_n}\right).\label{eq:13}
\end{eqnarray}
Let us first state the following technical lemma. 

\begin{lemma}\label{L:1} 
Let $\left(v_n\right)\in \mathcal{GS}\left(v^*\right)$, 
$\left(\beta_n\right)\in \mathcal{GS}\left(-\beta\right)$, and $m>0$ 
such that $m-v^*\xi>0$ where $\xi$ is defined in~\eqref{eq:xi}. We have 
\begin{eqnarray}\label{eq:14}
\lim_{n \to +\infty}v_nQ_n^{m}\sum_{k=1}^nQ_k^{-m}\frac{\beta_k}{v_k}
=\frac{1}{m-v^*\xi}. 
\end{eqnarray}
Moreover, for all positive sequence $\left(\alpha_n\right)$ such that $\lim_{n \to +\infty}\alpha_n=0$, and all $\delta \in \mathbb{R}$,
\begin{eqnarray}\label{eq:15}
\lim_{n \to +\infty}v_nQ_n^{m}\left[\sum_{k=1}^n Q_k^{-m} \frac{\beta_k}{v_k}\alpha_k+\delta\right]=0.
\end{eqnarray}
\end{lemma}

Lemma \ref{L:1} is widely applied throughout the proofs. Let us underline that it is its application, which requires Assumption $(A2)(iii)$ on the limit of $(n\gamma_n)$ as $n$ goes to infinity.\\

Our proofs are organized as follows. Propositions \ref{prop:bias:var:rn} and \ref{prop:MISE:rn} in Sections \ref{Section:prop:bias:var:rn} 
and \ref{Section:prop:MISE:rn} respectively, Theorem \ref{theo:TLC1} in Section \ref{Section:theo:TLC}. Propositions \ref{prop:bias:var:rn:bis} and \ref{prop:MISE:var:rn} in Sections \ref{Section:prop:bias:var:rn:bis} 
and \ref{Section:prop:MISE:rn:bis} respectively, Theorem \ref{theo:TLC2} in Section \ref{Section:theo:TLC:bis}.

\subsection{Proof of Proposition~\ref{prop:bias:var:rn}} \label{Section:prop:bias:var:rn}
Let us first note that, for $x$ such that $f_n\left(x\right)\not=0$, we have
\begin{eqnarray}\label{eq:decomp}
r_n\left(x\right)-r\left(x\right)&=&B_n\left(x\right)\frac{f\left(x\right)}{f_n\left(x\right)},
\end{eqnarray}
with 
\begin{eqnarray}\label{eq:Bn}
B_n\left(x\right)&=&\frac{1}{f\left(x\right)}\left(a_n\left(x\right)-a\left(x\right)\right)-\frac{r\left(x\right)}{f\left(x\right)}\left(f_n\left(x\right)-f\left(x\right)\right).
\end{eqnarray}
It follows from~(\ref{eq:decomp}), that the asymptotic behaviour of $r_n\left(x\right)-r\left(x\right)$ can be deduced from the one of  $B_n\left(x\right)$. Moreover, the following Lemma follows from the Proposition 1 of~\citet{Mok09a}.
\begin{lemma}[Bias and variance of $f_n$]
Let Assumptions $\left(A1\right)-\left(A3\right)$ and suppose that the stepsize $\left(\gamma_n\right)=\left(n^{-1}\right)$. 
\begin{enumerate}
\item If $a\in ]0, 1/5]$, then
\begin{eqnarray}\label{bias:fn}
\mathbb{E}\left[f_n\left(x\right)\right]-f\left(x\right)=\frac{1}{2\left(1-2a\right)}f^{\left(2\right)}\left(x\right)h_n^2\mu_2\left(K\right)+o\left(h_n^2\right).
\end{eqnarray}
If $a\in  ]1/5, 1[$, then
\begin{eqnarray}\label{bias:fn:bis}
\mathbb{E}\left[f_n\left(x\right)\right]-f\left(x\right)=o\left(\sqrt{n^{-1}h_n^{-1}}\right).
\end{eqnarray}
\item If $a\in [1/5, 1[$, then
\begin{eqnarray}
Var\left[f_n\left(x\right)\right]&=&\frac{1}{1+a}\frac{1}{nh_n}f\left(x\right)R\left(K\right)+o\left(\frac{1}{nh_n}\right).\label{var:fn}
\end{eqnarray}
If $a\in ]0,1/5[$, then
\begin{eqnarray}\label{var:fn:rep}
Var\left[f_n\left(x\right)\right]=o\left(h_n^4\right).
\end{eqnarray}
\end{enumerate}
\end{lemma}
Following similar steps as the proof of the Proposition 1 of~\citet{Mok09a}, we show that
\begin{lemma}[Bias and variance of $a_n$]\label{prop:bias:var:an}
Let Assumptions $\left(A1\right)-\left(A3\right)$ hold.
\begin{enumerate}
\item If $a\in ]0, \beta/5]$, then
\begin{eqnarray}\label{bias:an}
\mathbb{E}\left[a_n\left(x\right)\right]-a\left(x\right)=\frac{1}{2\left(1-2a\xi\right)}a^{\left(2\right)}\left(x\right)h_n^2\mu_2\left(K\right)+o\left(h_n^2\right).
\end{eqnarray}
If $a\in  ]\beta/5, 1[$, then
\begin{eqnarray}\label{bias:an:bis}
\mathbb{E}\left[a_n\left(x\right)\right]-a\left(x\right)=o\left(\sqrt{\beta_nh_n^{-1}}\right).
\end{eqnarray}
\item If $a\in [\beta/5, 1[$, then
\begin{eqnarray}
Var\left[a_n\left(x\right)\right]&=&\frac{\mathbb{E}\left[Y^2\vert X=x\right]f\left(x\right)}{\left(2-\left(\beta-a\right)\xi\right)}\frac{\beta_n}{h_n}R\left(K\right)+o\left(\frac{\beta_n}{h_n}\right).\label{var:an}
\end{eqnarray}
If $a\in ]0,\beta/5[$, then
\begin{eqnarray}\label{var:an:rep}
Var\left[a_n\left(x\right)\right]=o\left(h_n^4\right).
\end{eqnarray}
\end{enumerate}
\end{lemma} 
Then,~(\ref{bias:rn}) follows from~(\ref{bias:fn}),~(\ref{bias:an}) and~(\ref{eq:decomp}) and~(\ref{bias:rn:bis}) follows from~(\ref{bias:fn:bis}),~(\ref{bias:an:bis}) and~(\ref{eq:decomp}).\\
Now, it follows from~(\ref{eq:Bn}) that
\begin{eqnarray}\label{var:Bn}
\lefteqn{Var\left[B_n\left(x\right)\right]}\nonumber\\
&=&\frac{1}{f^2\left(x\right)}
\left\{Var\left[a_n\left(x\right)\right]+r^2\left(x\right)Var\left[f_n\left(x\right)\right]-2r\left(x\right)Cov\left(a_n\left(x\right),f_n\left(x\right)\right)\right\}.
\end{eqnarray}
In view of $\left(A3\right)$, and with the choice of the stepsize $\left(\gamma_n\right)=\left(n^{-1}\right)$ and using Lemma~\ref{L:1}, classical computations gives
\begin{eqnarray}\label{cov:an:fn}
Cov\left(a_n\left(x\right),f_n\left(x\right)\right)&=&\frac{\xi}{1+a\xi}\frac{\beta_n}{h_n}r\left(x\right)f\left(x\right)R\left(K\right)+o\left(\frac{\beta_n}{h_n}\right).
\end{eqnarray}
Then, the combination of~(\ref{eq:decomp}),~(\ref{var:Bn}),~(\ref{var:fn}),~(\ref{var:an}) and~(\ref{cov:an:fn}), gives~(\ref{var:rn}), and the combination of~(\ref{eq:decomp}),~(\ref{var:Bn}),~(\ref{var:fn:rep}),~(\ref{var:an:rep}) and~(\ref{cov:an:fn}), gives~(\ref{var:rn:rep}).

\subsection{Proof of Proposition~\ref{prop:MISE:rn}} \label{Section:prop:MISE:rn}
Following similar steps as the proof of the Proposition 2 of~\citet{Mok09a}, we proof the Propostion~\ref{prop:MISE:rn}.

\subsection{Proof of Theorem~\ref{theo:TLC1}} \label{Section:theo:TLC}
Let us at first assume that, if $a\geq\beta/5$, then 
\begin{eqnarray}\label{eq:22}
\sqrt{\beta_n^{-1} h_n}\left(r_{n}\left( x\right)-\mathbb{E}\left[r_n\left(x\right)\right]\right)\stackrel{\mathcal{D}}{\rightarrow}\mathcal{N}\left( 0,
V_{a,\xi,\beta}^{\left(1\right)}\right).
\end{eqnarray}
In the case when $a>\beta/5$, Part 1 of Theorem~\ref{theo:TLC1} follows from the combination of~\eqref{bias:rn:bis} and~\eqref{eq:22}. In the case when $a=\beta/5$, Parts 1 and 2 of Theorem~\ref{theo:TLC1} follow from the combination of~\eqref{bias:rn} and~\eqref{eq:22}. In the case $a<\beta/5$,~\eqref{var:rn:rep} implies that 
\begin{eqnarray*}
h_n^{-2}\left(r_n\left(x\right)-\mathbb{E}\left(r_n\left(x\right)\right)\right)\stackrel{\mathbb{P}}{\rightarrow}0,
\end{eqnarray*}
and the application of~\eqref{bias:rn} gives Part 2 of Theorem~\ref{theo:TLC1}.\\

We now prove~\eqref{eq:22}. 
In view of~\eqref{eq:Bn}, we have
\begin{eqnarray}\label{eq:Bn:Ebn}
B_n\left(x\right)-\mathbb{E}\left[B_n\left(x\right)\right]
&=&\frac{1}{f\left(x\right)}Q_n\sum_{k=1}^n\left(T_k\left(x\right)-\mathbb{E}\left[T_k\left(x\right)\right]\right),
\end{eqnarray}
with 
\begin{eqnarray}\label{eq:Tk}
T_{k}\left(x\right)&=&Q_k^{-1}\left(\beta_kZ_k\left(x\right)-r\left(x\right)\zeta_n\zeta_k^{-1}\gamma_kW_k\left(x\right)\right).
\end{eqnarray}
In the case when $\left(\gamma_n\right)=\left(n^{-1}\right)$, we have $\zeta_n=\left(nQ_n\right)^{-1}$ et $\zeta_k^{-1}\gamma_k=Q_k$, then
\begin{eqnarray*}
T_{k}\left(x\right)&=&Q_k^{-1}\beta_kZ_k\left(x\right)-r\left(x\right)\left(nQ_n\right)^{-1}W_k\left(x\right).
\end{eqnarray*}
Set
\begin{eqnarray}\label{eq:Yk}
Y_k\left(x\right)&=&T_{k}\left(x\right)-\mathbb{E}\left(T_{k}\left(x\right)\right).
\end{eqnarray}
Moreover, we have
\begin{eqnarray*} 
v_n^2&=&\sum_{k=1}^nVar\left(Y_{k}\left(x\right)\right)\nonumber\\
&=&\sum_{k=1}^nQ_k^{-2}\beta_k^2Var\left(Z_k\left(x\right)\right)+r^2\left(x\right)\left(nQ_n\right)^{-2}\sum_{k=1}^nVar\left(W_k\left(x\right)\right)\nonumber\\
&&-2r\left(x\right)\left(nQ_n\right)^{-1}\sum_{k=1}^nQ_k^{-1}\beta_kCov\left(Z_k\left(x\right),W_k\left(x\right)\right).
\end{eqnarray*}
Moreover, in view of $\left(A3\right)$, classical computations give
\begin{eqnarray*}
Var\left(Z_k\left(x\right)\right)&=&\frac{1}{h_k}\left[\mathbb{E}\left[Y^2\vert X=x\right]f\left(x\right)R\left(K\right)+o\left(1\right)\right],\\
Var\left(W_k\left(x\right)\right)&=&\frac{1}{h_k}\left[f\left(x\right)R\left(K\right)+o\left(1\right)\right],\\
Cov\left(Z_k\left(x\right),W_k\left(x\right)\right)&=&\frac{1}{h_k}\left[r\left(x\right)f\left(x\right)R\left(K\right)+o\left(1\right)\right].
\end{eqnarray*}
The application of Lemma $\ref{L:1}$ ensures that
\begin{eqnarray*}
v_n^2&=&\sum_{k=1}^n\frac{Q_k^{-2}\beta_k^2}{h_k}\left[\mathbb{E}\left[Y^2\vert X=x\right]f\left(x\right)R\left(K\right)+o\left(1\right)\right]\nonumber\\
&&+\frac{r\left(x\right)}{n^2Q_n^2}\sum_{k=1}^n\frac{1}{h_k}\left[f\left(x\right)R\left(K\right)+o\left(1\right)\right]\nonumber\\
&&-2\frac{r\left(x\right)}{nQ_n}\sum_{k=1}^n\frac{Q_k^{-1}\beta_k}{h_k}\left[r\left(x\right)f\left(x\right)R\left(K\right)+o\left(1\right)\right]\\
&=&\frac{f^2\left(x\right)}{Q_n^2}\frac{\beta_n}{h_n}\left[V_{a,\xi,\beta}^{\left(1\right)}+o\left(1\right)\right].
\end{eqnarray*}
On the other hand, we have, for all $p>0$, 
\begin{eqnarray*}
\mathbb{E}\left[\left|T_k\left(x\right)\right|^{2+p}\right] &=&
O\left(\frac{1}{h_k^{1+p}}\right),
\end{eqnarray*}
and, since $\lim_{n\to\infty}\left(n\beta_n\right)>\left(\beta-a\right)/2$, there exists $p>0$ such that $\lim_{n\to \infty}\left(n\beta_n\right)>\frac{1+p}{2+p}\left(\beta-a\right)$. Applying Lemma $\ref{L:1}$, we get 
\begin{eqnarray*}
\sum_{k=1}^n\mathbb{E}\left[\left|Y_{k}\left(x\right)\right|^{2+p}\right]&=&O\left(\sum_{k=1}^n Q_k^{-2-p}\beta_k^{2+p}\mathbb{E}\left[\left|T_k\left(x\right)\right|^{2+p}\right]\right)\nonumber\\
&=&O\left(\sum_{k=1}^n \frac{Q_k^{-2-p}\beta_k^{2+p}}{h_k^{1+p}}\right)\\
&=&O\left(\frac{\beta_n^{1+p}}{Q_n^{2+p}h_n^{1+p}}\right)\nonumber,
\end{eqnarray*}
and we thus obtain 
\begin{eqnarray*}
\frac{1}{v_n^{2+p}}\sum_{k=1}^n\mathbb{E}\left[\left|Y_{k}\left(x\right)\right|^{2+p}\right]& = & O\left({\left[\beta_nh_n^{-1}\right]}^{p/2}\right)=o\left(1\right).
\end{eqnarray*}
The convergence in~\eqref{eq:22} then follows from the application of Lyapounov's Theorem.

\subsection{Proof of Proposition~\ref{prop:bias:var:rn:bis}}\label{Section:prop:bias:var:rn:bis}
The following Lemma follows from the Proposition 1 of~\citet{Mok09a}.
\begin{lemma}[Bias and variance of $f_n$]
Let Assumptions $\left(A1\right)-\left(A3\right)$ and suppose that the stepsize $\left(\gamma_n\right)=\left(\left[1-a\right]n^{-1}\right)$. 
\begin{enumerate}
\item If $a\in ]0, 1/5]$, then
\begin{eqnarray}\label{bias:fn:2}
\mathbb{E}\left[f_n\left(x\right)\right]-f\left(x\right)=\frac{1-a}{2\left(1-3a\right)}f^{\left(2\right)}\left(x\right)h_n^2\mu_2\left(K\right)+o\left(h_n^2\right).
\end{eqnarray}
If $a\in  ]1/5, 1[$, then
\begin{eqnarray}\label{bias:fn:bis:2}
\mathbb{E}\left[f_n\left(x\right)\right]-f\left(x\right)=o\left(\sqrt{n^{-1}h_n^{-1}}\right).
\end{eqnarray}
\item If $a\in [1/5, 1[$, then
\begin{eqnarray}
Var\left[f_n\left(x\right)\right]&=&\frac{1-a}{nh_n}R\left(K\right)+o\left(\frac{1}{nh_n}\right).\label{var:fn:2}
\end{eqnarray}
If $a\in ]0,1/5[$, then
\begin{eqnarray}\label{var:fn:rep:2}
Var\left[f_n\left(x\right)\right]=o\left(h_n^4\right).
\end{eqnarray}
\end{enumerate}
\end{lemma}
Then,~(\ref{bias:rn:min:var}) follows from~(\ref{bias:fn:2}),~(\ref{bias:an}) and~(\ref{eq:decomp}) and~(\ref{bias:rn:bis:min:var}) follows from~(\ref{bias:fn:bis:2}),~(\ref{bias:an:bis}) and~(\ref{eq:decomp}).\\
Moreover, in view of $\left(A3\right)$, and using the choice of the stepsize $\left(\gamma_n\right)=\left(\left[1-a\right]n^{-1}\right)$ and using Lemma~\ref{L:1}, classical computations gives
\begin{eqnarray}\label{cov:an:fn:2}
Cov\left(a_n\left(x\right),f_n\left(x\right)\right)&=&\left(1-a\right)\xi\frac{\beta_n}{h_n}r\left(x\right)f\left(x\right)R\left(K\right)+o\left(\frac{\beta_n}{h_n}\right).
\end{eqnarray}
Then, the combination of~(\ref{eq:decomp}),~(\ref{var:Bn}),~(\ref{var:fn}),~(\ref{var:an}) and~(\ref{cov:an:fn:2}), gives~(\ref{var:rn:min:var}), and the combination of~(\ref{eq:decomp}),~(\ref{var:Bn}),~(\ref{var:fn:rep:2}),~(\ref{var:an:rep}) and~(\ref{cov:an:fn:2}), gives~(\ref{var:rn:rep:min:var}).

\subsection{Proof of Proposition~\ref{prop:MISE:var:rn}} \label{Section:prop:MISE:rn:bis}
Following similar steps as the proof of the Proposition 2 of~\citet{Mok09a}, we proof the Propostion~\ref{prop:MISE:var:rn}.
\subsection{Proof of Theorem~\ref{theo:TLC2}} \label{Section:theo:TLC:bis}
Following similar steps as the proof of the Theorem~\ref{theo:TLC1} and using the fact that in the case when $\left(\gamma_n\right)=\left(\left[1-a\right]n^{-1}\right)$, we have $Q_k^{-1}\zeta_n\zeta_k^{-1}\gamma_k=\left(1-a\right)h_k/(nh_nQ_n)$, and then it follows from~(\ref{eq:Tk}), that
\begin{eqnarray*}
T_{k}\left(x\right)&=&Q_k^{-1}\beta_kZ_k\left(x\right)-r\left(x\right)\frac{\left(1-a\right)h_k}{nh_nQ_n}W_k\left(x\right),
\end{eqnarray*}
we prove Theorem~\ref{theo:TLC2}.

\section*{}
\makeatletter
\renewcommand{\@biblabel}[1]{}
\makeatother

$ $\\
Universit\'e de Poitiers, Laboratoire de Math\'ematiques et Application, Futuroscope Chasseneuil, France\\
E.mail: Yousri.Slaoui@math.univ-poitiers.fr

\begin{thebibliography}{99}

\bibitem[{Altman and Leger(1995)}]{Alt95}
{Altman, N.} and {Leger, C.} (1995).
\newblock {Bandwidth selection for kernel distribution function estimation}.
\newblock \textit{J. Statist. Plann. Inference}, {\bf 46}, 195--214.


\bibitem[{Bojanic and Seneta(1995)}]{Boj73}
{Bojanic, R.} and {Seneta, E.} (1973).
\newblock {A unified theory of regularly varying sequences}.
\newblock \textit{Math. Z.}, {\bf 134}, 91--106.

\bibitem[{Delaigle and Gijbels(2004)}]{Del04}
{Delaigle, A.} and {Gijbels, I.} (2004). 
\newblock {Practical bandwidth selection in deconvolution kernel density estimation}.
\newblock \textit{Comput. Statist. Data Anal.}, {\bf 45}, 249--267.

\bibitem[{Galambos and Seneta(1973)}]{Gal73}
{Galambos, J.} and {Seneta, E.} (1973).
\newblock {Regularly varying sequences}.
\newblock \textit{Proc. Amer. Math. Soc.}, {\bf 41}, 110--116.

\bibitem[{Hart and Vieu(1990)}]{Har90}
{Hart, J. D.} and {Vieu, P.} (1990).
\newblock {Data-Driven Bandwidth Choice for Density Estimation Based
on Dependent Data}. 
\newblock \textit{Ann. Statist.}, {\bf 18}, 873--890.

\bibitem[{Huang et al.(2014)}]{Hua14}
{Huang, Y.}, {Chen, X.} and {Wu, W.~B.} (2014).
\newblock {Recursive nonparametric estimation for times series}.
\newblock \textit{IEEE Trans. Inform. Theory}, {\bf 60}, 1301--1312.

\bibitem[{Mokkadem and Pelletier(2007)}]{Mok07}
{Mokkadem, A.} and {Pelletier, M.} (2007).
\newblock {A companion for the Kiefer-Wolfowitz-Blum stochastic approximation algorithm}.
\newblock \textit{Ann. Statist.}, {\bf 35}, 1749--1772.


\bibitem[{Mokkadem et~al.(2009a)Mokkadem, Pelletier and Slaoui(2009)}]{Mok09a} 
{Mokkadem, A.} {Pelletier, M.} and {Slaoui, Y.} (2009a).
\newblock {The stochastic approximation method for the estimation of a multivariate probability density}.
\newblock \textit{J. Statist. Plann. Inference}, {\bf 139}, 2459--2478.

\bibitem[{Mokkadem et~al.(2009b)Mokkadem, Pelletier and Slaoui(2009)}]{Mok09b} 
{Mokkadem, A.} {Pelletier, M.} and {Slaoui, Y.} (2009b).
\newblock {Revisiting R\'ev\'esz's stochastic approximation method for the estimation of a regression function}.
\newblock \textit{ALEA. Lat. Am. J. Probab. Math. Stat.}, {\bf 6}, 63--114.


\bibitem[Nadaraya(1964)]{Nad64} 
{Nadaraya, E. A.} (1964).
\newblock {On estimating regression}.
\newblock \textit{Theory Probab. Appl.}, {\bf 10}, 186--190.

\bibitem[R\'ev\'esz(1973)]{Rev73}
{R\'ev\'esz, P.} (1973).
\newblock {Robbins-Monro procedure in a Hilbert space and its
application in the theory of learning processes I}.
\newblock \textit{Studia Sci. Math. Hung.}, {\bf 8}, 391--398.

\bibitem[R\'ev\'esz(1977)]{Rev77}
{R\'ev\'esz, P.} (1977).
\newblock {How to apply the method of stochastic approximation in the non-parametric estimation of a regression function}.
\newblock \textit{Math. Operationsforsch. Statist., Ser. Statistics}, {\bf 8}, 119--126.

\bibitem[{Silverman(1986)}]{Sil86} 
{Silverman, B.~W.} (1986).
\newblock {Density estimation for statistics and data analysis}.
\newblock \textit{Chapman and Hall}, London.


\bibitem[{Slaoui(2013)}]{Sla13} 
{Slaoui, Y.} (2013).
\newblock {Large and moderate deviation principles for recursive kernel density estimators defined by stochastic approximation method}.
\newblock \textit{Serdica Math. J.} {\bf 39},~53--82.

\bibitem[{Slaoui(2014a)}]{Sla14a} 
{Slaoui, Y.} (2014a).
\newblock {Bandwidth selection for recursive kernel density estimators defined by stochastic approximation method}.
\newblock \textit{J. Probab. Stat}, {\bf 2014}, ID 739640, doi:10.1155/2014/739640.

\bibitem[{Slaoui(2014b)}]{Sla14b} 
{Slaoui, Y.} (2014b).
\newblock {The stochastic approximation method for the estimation of a distribution function}.
\newblock \textit{Math. Methods Statist.} {\bf 23},~306--325.

\bibitem[{Slaoui(2015a)}]{Sla15a} 
{Slaoui, Y.} (2015a).
\newblock {Moderate deviation principles for recursive regression estimators defined by stochastic approximation method}.
\newblock \textit{Int. J. Math. Stat.} {\bf 16},~51--60.

\bibitem[{Slaoui(2015b)}]{Sla15b} 
{Slaoui, Y.} (2015b).
\newblock {Large and moderate deviation principles for recursive regression estimators defined by stochastic approximation method}.
\newblock \textit{Serdica Math. J.} {\bf 41},~307--328.


\bibitem[{Slaoui(2015c)}]{Sla15c} 
{Slaoui, Y.} (2015c).
\newblock {Plug-In Bandwidth selector for recursive kernel regression estimators defined by stochastic approximation method.}
\newblock \textit{Stat. Neerl.} {\bf 69},~483--509.

\bibitem[{Tsybakov(1990)}]{Tsy90} 
{Tsybakov, A.~B.} (1990).
\newblock {Recurrent estimation of the mode of a multidimensional distribution}.
\newblock \textit{Probl. Inf. Transm.}, {\bf 8}, 119--126.

\bibitem[Watson(1964)]{Wat64} 
{Watson, G. S.} (1964). 
\newblock {Smooth regression analysis}. 
\newblock \textit{Sankhya A}, {\bf 26}, 359--372.

\end{thebibliography}
\end{document}